\documentclass[10pt,reqno]{amsart}
\newtheorem{MainThm}{Theorem}

\usepackage{amssymb, amsmath, euscript, verbatim, sw, default, array, url, mathrsfs}
\usepackage[margin=32pt,font=footnotesize]{caption}

\begin{document}
  \title[Lagrangian Correspondences and Donaldson's TQFT]{Lagrangian Correspondences and Donaldson's TQFT Construction of the Seiberg-Witten Invariants of $3$-Manifolds}
  \author{Timothy Nguyen}
  \date{\today}

\begin{abstract}
  Using Morse-Bott techniques adapted to the gauge-theoretic setting, we show that the limiting boundary values of the space of finite energy monopoles on a connected $3$-manifold with at least two cylindrical ends provides an immersed Lagrangian submanifold of the vortex moduli space at infinity.  By studying the signed intersections of such Lagrangians, we supply the analytic details of Donaldson's TQFT construction of the Seiberg-Witten invariants of a closed $3$-manifold.
\end{abstract}

\maketitle

\tableofcontents

\section{Introduction}

In \cite{Don99}, Donaldson outlines a 2+1 topological quantum field theoretic construction of the Seiberg-Witten invariants of closed (oriented) $3$-manifolds with $b_1 > 0$.  This construction is motivated from the expectation that the moduli space of finite energy monopoles on a $3$-manifold with cylindrical ends provides an (immersed) Lagrangian correspondence between the vortex moduli spaces at infinity.  From this expectation, by regarding the closed $3$-manifold $Y$ as a cobordism $W$ with its two boundary components identified, and then decomposing $W$ into a composite of cobordisms, one can interpret the Seiberg-Witten invariants of $Y$ as being obtained from the composite of the cyclic sequence of Lagrangian correspondences.

In this paper, we supply the analytic details of these expectations.  Our first main result is the following theorem which establishes the previously described Lagrangian correspondence under suitable assumptions:

\begin{MainThm}
  Let $W$ be a connected oriented $3$-manifold with at least two cylindrical ends. Then for any $\spinc$ structure and suitable generic coclosed perturbations of the Seiberg-Witten equations, the moduli space of finite energy monopoles on $W$ is a smooth, compact, orientable manifold. Moreover, the map which sends a monopole to its limiting value at infinity along the ends sends the moduli space of monopoles to an immersed Lagrangian inside the vortex moduli spaces at infinity.
\end{MainThm}

A more precise formulation, which takes into account the nature of the coclosed perturbations, is to be found in Theorem \ref{ThmOne}. The assumption about the underlying three-manifold having at least two boundary components is to ensure that we can perturb the Seiberg-Witten equations in such a way that at the limiting ends the equations can be interpreted as a Morse-Bott flow. (Equivalently, we can perturb so that the limiting vortex moduli space contains no reducibles.) To the author's knowledge, Theorem 1 (or some version thereof) has essentially been a folk theorem since the early days of Seiberg-Witten theory.  Indeed, the conjectured equivalence between Heegaard Floer homology and Seiberg-Witten Floer homology, recently established in \cite{KLT}, is credible if one expects that the space of monopoles on the handlebodies occurring in a Heegaard decomposition of a closed $3$-manifold yield Lagrangians in the vortex moduli space on the Heegaard surface (in the limit in which the neck surrounding the surface is stretched to infinity). This is only a heuristic picture of course, since the vortex moduli space one obtains from the Seiberg-Witten setup from a genus $g$ handlebody is a symmetric product of order $g - 1$ (see Lemma \ref{LemmaVortices}) and not of order $g$ as would be needed for Heegaard-Floer theory.


Our second main theorem, makes precise the invariant that Donaldson's construction in \cite{Don99} computes. First some notation. Let $Y$ be a closed oriented $3$-manifold with $b_1(Y) > 0$.  Pick any connected nonseparating orientable hypersurface $\Sigma \subset Y$ and form the cylindrical end manifold $W^*$ from the manifold $W = Y\setminus\Sigma$ by attaching two semi-infinite ends $(-\infty,0] \times -\Sigma$ and $[0,\infty) \times \Sigma$ in the obvious way. Pick a $\spinc$ structure $\s_0$ on $W^*$. The metric and $\spinc$ structure on $W^*$ are assumed to be a product in the natural way on the semi-infinite ends. Regard $Y$ as $W$ with its two ends identified by a diffeomorphism $h: \Sigma \to \Sigma$, and let $\Spin^c(Y,\s_0)$ denote the set of all $\spinc$ structures on $Y$ obtained from $\s_0|_W$ by all possible ways of gluing along $\Sigma$.

\begin{MainThm}
  Let $Y$ be a closed oriented $3$-manifold with $b_1(Y) > 0$.  For suitably chosen generic perturbations $\eta \in \Omega^1(Y; i\R)$, let $[L]$ denote the homology class of the Lagrangian obtained from Theorem 1 applied to $W^*$. Then we have
  \begin{equation}
    \sum_{\s \in \Spinc(Y,\s_0)} SW(\s, \eta) = [L] \cap [\Gamma_h] \label{intersect0}
  \end{equation}
  where the left-hand side is a sum of the Seiberg-Witten invariants of $Y$ with respect to the perturbation $\eta$ and $\spinc$ structures belonging to $\Spin(Y,\s_0)$, and the right-hand side denotes the signed intersection of $[L]$ and the graph $\Gamma_h$ induced by $h$ inside the vortex moduli space on $-\Sigma \times \Sigma$. Here, a homology orientation on $Y$ and an orientation of $[L]$ are chosen compatibly (each of these determines an overall sign for the left-hand side and right-hand side, respectively).
\end{MainThm}

See Theorem \ref{ThmTwo} for a more precise formulation. Formula (\ref{intersect0}) is precisely that which appears in \cite{Mark} without proof (although see Remark \ref{RemMark}), which T. Mark uses to prove a (version of a) conjecture of Hutchings and Lee \cite{HL}.  We should remark that while this formula at face value appears difficult to explicitly compute, its important feature is that it is obtained from a 2+1 TQFT-like framework.  That is, there is an underlying composition rule, whereby if we can decompose the $3$-manifold $Y$ into simple pieces, namely, a composite of elementary cobordisms with the incoming and outgoing ends identified, our formula (\ref{intersect}) is obtained from understanding the composite of the morphisms produced from the elementary cobordisms.  Donaldson, using only algebraic and topological arguments, provides an explicit and elegant computation for these maps induced from elementary cobordisms which T. Mark then exploits in \cite{Mark}.  Moreover, Donaldson himself uses his computations to recover the formula of Meng-Taubes \cite{MT}, relating the Seiberg-Witten invariants to the Alexander polynomial, in the case of $b_1 = 1$. In this way, we expect formula (\ref{intersect}) to be a useful addition to the list of ways one can compute and interpret the Seiberg-Witten invariants of $3$-manifolds (see also \cite{Tur}). In fact, Theorem \ref{ThmTwo} proves more than just formula (\ref{intersect}), see Remark \ref{RemGen}.

We conclude this introduction with a summary of the ideas involved in the proofs.  For Theorem \ref{ThmOne}, we proceed by first analyzing the (perturbed) Seiberg-Witten equations on a semi-infinite cylinder $[0,\infty) \times \Sigma$, which one can interpret (formally) as a Morse-Bott flow restricted to the level set of a moment map.  From this, we adapt standard Morse-Bott techniques, which appear to be well-documented in the instanton literature (\cite{MMR}, \cite{Don99}), to the Seiberg-Witten case.  This allows us to give an explicit description of the moduli space of finite energy monopoles on a semi-infinite cylinder, namely, that it is the ``stable manifold" to the space of vortices on $\Sigma$ under the Morse-Bott flow. We then piece together the moduli space of monopoles arising from semi-infinite cylinders and from compact $3$-manifolds with boundary (see \cite{N1}) via a suitable fiber product to obtain the moduli space of monopoles on a general $3$-manifold with cylindrical ends. The symplectic properties of these moduli spaces we obtain depend upon fundamental properties of Dirac operators and their boundary values, where a general technique known as the ``invertible double" is used in the context of weighted spaces.  We also establish as easy consequences Theorems \ref{ThmSmallEnergy} and \ref{ThmLagVort}, which provide Lagrangian submanifolds in the space of connections and spinors on $\Sigma$ whose topological type is explicit.  Applications of this include the study of Lagrangian boundary conditions \cite{N2} and semi-infinite cycles \cite{Lip}.

For Theorem \ref{ThmTwo}, there are two main technical ingredients involved.  The first is the issue of gluing moduli spaces in the Morse-Bott setting.  We should note that because we are gluing together cobordisms, the result of which is another cobordism with possibly nonempty boundary, some care must be taken since it is incorrect to pass to the cobordism times $S^1$ and use four-dimensional gluing.  (Indeed, a monopole on a non-closed $3$-manifold times $S^1$ need not be pulled back from a monopole on the $3$-manifold.)  However, by the proper use of weighted spaces and some care with asymptotic boundary conditions, one can proceed with gluing along lines similar to the Morse nondegenerate case. The second and main technical issue however is the issue with signs in formula (\ref{intersect}), namely, why the ``geometric signs", arising from the signed intersection of Lagrangian correspondences, agree with the ``analytic signs", arising from orienting determinant lines of Fredholm operators that occur in counting monopoles.  Here, one needs to understand how to orient the moduli space of monopoles on a cylindrical end manifold and how to glue these orientations together in the Morse-Bott framework.  The cylindrical end nature along with the Morse-Bott situation makes the orientation issue delicate, since there is nontrivial topology at infinity.  Moreover, our configuration space consists of configurations asymptotic to vortices that are not reducible, which makes orientating determinant lines not straightforward.  Nevertheless, an excision argument combined with the fact that the gauge-fixed linearized vortex operator is complex linear with respect to a suitably chosen complex structure allows us to establish orientability of the monopole moduli space.

Let us remark that there is another way in which one could analyze the issue of signs.  Indeed, in the instanton case, there is Taubes's gauge theoretic construction of Casson's invariant \cite{T}, which equates the signed intersection of Lagrangian submanifolds in the space of flat connections on a Heegaard surface to a signed count of flat connections on the relevant $3$-manifold.  The essential ingredient in Taubes's work is to equate spectral flow (analytic) with a particular Maslov index (geometric) modulo two.  The equivalence between these two quantities in the general integer setting was later recast in greater generality for (neck cylindrical) Dirac operators by Nicolaescu \cite{Nic95}. If we were to try to adapt the approach of Taubes directly to the Seiberg-Witten setting,  we would a priori have to deal with the Maslov index in infinite dimensions, since the boundary values of the space of monopoles yields infinite dimensional Lagrangians.  However, since the Maslov index is preserved under symplectic reduction under suitable hypotheses, we can relate the infinite dimensional Maslov index (obtained from studying monopoles on elementary pieces of the closed $3$-manifold) to the finite dimensional Maslov index (stretching the neck of these elementary pieces to infinity) occurring in the vortex moduli spaces.

In fact, the author first went about proving Theorem \ref{ThmTwo} in the above way, combining Morse-Bott techniques with a generalization of Nicolaescu's results to non-cylindrical Dirac operators.  However, the resulting analysis becomes more technical than the one presented here.  Moreover, the main drawback concerning the use of spectral flow techniques is that it is limited to self-adjoint operators whereas the gluing and orientation methods presented here only require that we are in the more general Morse-Bott setting.  Indeed, our proof here was motivated by the fact that in \cite{MMS} a similar but unproven claim appears about the equivalence of the 4-dimensional Seiberg-Witten invariant and a signed intersection of manifolds at infinity, see \cite[Theorem 4.1]{MMS}.  In fact, one can interpret most of the literature on the gluing of Seiberg-Witten invariants as variations on the same theme of signed intersections, albeit the only cases we could find in the literature supplied with a proof are ones in which the critical sets at infinity are Morse nondegenerate (and hence the signed intersection reduces to a signed product formula).  The convenience of working in the present three-dimensional setting is that the boundaries of three-manifolds, being two-manifolds, are completely understood along with their corresponding critical sets for the Chern-Simons-Dirac functional, the vortex moduli spaces.  It is only because of this explicit description that we are able to establish orientability of the moduli space of monopoles on cylindrical end $3$-manifolds.  In principle, one could adapt the methods here to prove signed intersection formulas for the four-dimensional Seiberg-Witten invariants in the Morse-Bott situation (thus, generalizing the Morse nondegenerate situation of, e.g.,  \cite{MST}). However, because Morse-Bott critical sets for the Chern-Simons-Dirac functional are not explicitly computable for a general $3$-manifold, it is not a priori guaranteed that the moduli space of monopoles on a cylindrical end $4$-manifold is orientable (except in the standard Morse nondegenerate case).\\



\noindent \textit{Acknowledgements.} The author thanks Tom Mrowka, Liviu Nicolaescu, and Tim Perutz for providing valuable discussions. Perutz also deserves special recognition for having suggested and inspired this line of research.


\section{The Seiberg-Witten Invariant}

We quickly review the construction of the Seiberg-Witten invariant of a closed Riemannian $3$-manifold $Y$ with $b_1 > 0$. Given a $\spinc$ structure $\s$ on $Y$, we obtain the spinor bundle $\S = \S(\s)$ associated to $\s$, uniquely determined up to isomorphism by requiring that Clifford multiplication $\rho: \Lambda^*(TY) \to \mathrm{End}(\S)$ maps the volume form on $Y$ to the identity automorphism on $\S$.  Let $L = \det(\s)$ denote the determinant line bundle of $\S$. From this, we obtain the configuration space $\fC(Y) = \A(Y) \times \Gamma(\S)$, where $\A(Y)$ denotes the space of compatible $\spinc$ connections on $\S$ and $\Gamma(\S)$ is the space of smooth sections of $\S$.  We write $(B,\Psi)$ to denote the pair of a $\spinc$ connection and spinor on $Y$.

The Seiberg-Witten equations on $Y$ are given by
\begin{equation}
\begin{split}
  *F_B + \rho^{-1}(\Psi\Psi^*)_0 &= \eta \\
  D_B\Psi & = 0.
\end{split}\label{eq:SW3}
\end{equation}
In the first line, $F_B$ denotes half the curvature of the connection on the determinant line bundle $L$ induced by $B$. (If $L$ has a square root $L^{1/2}$, then $F_B$ would be the curvature of the connection induced on $L^{1/2}$).
Next, $*$ is the Hodge star operator on $Y$, and $(\Psi\Psi^*)_0$ is the traceless part of the Hermitian endomorphism $\Psi \otimes \Psi^* \in \mathrm{End}(\S)$.  The term $\eta \in \Omega^1(Y;i\R)$ is an imaginary coclosed form serving as a perturbation.  Finally, $D_B: \Gamma(\S) \to \Gamma(\S)$ is the $\spinc$ Dirac operator coupled to $B$.

By the standard theory, when $b_1(Y) > 0$, then for any metric on $Y$ and generic choice of $\eta$, the moduli space of gauge equivalence classes of solutions to (\ref{eq:SW3}) is a compact zero-dimensional moduli space.  An appropriate signed count of these solutions, which involves orienting determinant lines of families of Fredholm operators (see e.g. \cite{Nic}), gives us an integer which is independent of the choice of metric and generic perturbation when $b_1(Y) > 1$ and exhibits a wall-crossing phenomenon when $b_1(Y)=1$.  Thus, the Seiberg-Witten equations yield for us a map
$$SW: \Spinc(Y) \to \Z$$
assigning to each $\spinc$ structure on $Y$ the associated signed count of monopoles. This is well-defined when $b_1(Y) > 1$ and depends on a choice of chamber in $H^1(Y)$ when $b_1(Y) =1$.  The values that $SW$ produces as we vary the $\spinc$ structure on $Y$ (and the choice of chamber for $\eta$ when $b_1(Y)=1$) are known as the Seiberg-Witten invariants of $Y$.


\section{The Seiberg-Witten Flow}

Let $Y = [0,\infty) \times \Sigma$, where $\Sigma$ is a connected Riemann surface, and endow $Y$ with the product metric.  In this cylindrical situation, we will interpret the Seiberg-Witten equations, in temporal gauge, as a downward gradient flow of a Chern-Simons-Dirac functional on $\Sigma$ restricted to the level set of a moment map. Let $t \in  [0,\infty)$ be the time-variable.  Let $(B,\Psi)$ be a smooth solution to $SW_3(B,\Psi) = 0$ on $Y$.   We will always take the $\spinc$ structure on $Y$ to be pulled back from a $\spinc$ structure on $\Sigma$, and by abuse of notation, we denote both of these $\spinc$ structures by $\s$. With respect to this product structure, we can write the equations $SW_3(B,\Psi) = 0$ in a rather explicit fashion.  Recall that every Kahler manifold has a canonical $\spinc$ structure.  For a Riemann surface $\Sigma$, the spinor bundle associated to this canonical $\spinc$ structure is isomorphic to $K_\Sigma^{1/2} \oplus K_\Sigma^{-1/2}$, where $K_\Sigma$ is the canonical bundle of $\Sigma$.  Moreover, a $\spinc$ structure on $\Sigma$ is uniquely determined by its determinant line bundle $L$, and the corresponding spinor bundle it determines is isomorphic to
\begin{equation}
  \S_\Sigma \cong (K_\Sigma \otimes L)^{1/2} \oplus (K_\Sigma^{-1} \otimes L)^{1/2}. \label{eq2:Siso}
\end{equation}
Let $\pi_\Sigma: [0,\infty) \times \Sigma \to \Sigma$ denote the natural projection.  From the above, given a $\spinc$ structure on $ [0,\infty) \times \Sigma$ pulled back  from a $\spinc$ structure on $\Sigma$ via $\pi_\Sigma$, the spinor bundle $\S$ on $ [0,\infty) \times \Sigma$ is isomorphic to
\begin{equation}
  \pi_\Sigma^*(K_\Sigma \otimes L)^{1/2} \oplus \pi_\Sigma^*(K_\Sigma^{-1} \otimes L)^{1/2}. \label{eq2:spinor-decomp}
\end{equation}
From now on, we always assume that we are in this product situation on $[0,\infty) \times \Sigma$.  Since $T^*( [0,\infty) \times \Sigma) \cong T^*[0,\infty) \oplus T^*\Sigma$, we can always choose our Clifford multiplication $\rho$ on $ [0,\infty) \times \Sigma$ to be such that $\rho$ factors through the direct sum decomposition of $T^*([0,\infty) \times \Sigma)$.  From this, we can choose $\rho$ so that
$$\rho(\partial_t) \equiv \begin{pmatrix}i & 0 \\ 0 & -i\end{pmatrix}$$
with respect to the decomposition (\ref{eq2:spinor-decomp}). Using local holomorphic coordinates $z = x + iy$ on $\Sigma$, we can decompose a $1$-form on $Y$ into its $dt$, $dz$, and $d\bar z$ components.  Given a $\spinc$ connection $B$ on $Y$, let $F_{x,y}dx \wedge dy + F_{x,t}dx \wedge dt + F_{y,t}dy\wedge dt$ denote the local coordinate representation of $F_B$.  Then the equation $*F_B + \rho^{-1}(\Psi\Psi^*_0) = 0$ appearing in the unperturbed Seiberg-Witten equation $SW_3(B,\Psi) = 0$ can be written explicitly as (see \cite{MST})\footnote{Note that our sign conventions are that of \cite{KM}, namely $\rho(dtd\Sigma) = 1$, which is the opposite choice of sign in \cite{MST}. }:
\begin{align}
  \left(F_{x,y} + \frac{i}{2}(|\Psi_+|^2 - |\Psi_-|^2)\right)dt &= 0 \label{eq2:SW3-1}\\
  \frac{1}{2}(F_{y,t}-iF_{x,t})d\bar z + \bar\Psi_+\Psi_- & = 0 \label{eq2:SW3-2}\\
  \frac{1}{2}(F_{y,t}+iF_{x,t})dz + \Psi_+\bar\Psi_- & = 0. \label{eq2:SW3-3}
\end{align}
Here, $\Psi = (\Psi_+,\Psi_-)$ is the decomposition of $\Psi$ with respect to  (\ref{eq2:spinor-decomp}), so that $\bar\Psi_+\Psi_-$ and $\Psi_+\bar\Psi_-$ are well-defined elements of $\pi^*_\Sigma K_\Sigma^{\mp}$, respectively.  Observe that the last equation above is just the complex conjugate of the second.  Moreover, the Dirac equation $D_B\Psi = 0$ becomes
\begin{equation}
  \begin{pmatrix}
    i\nabla_{B,\partial_t} & \sqrt{2}\bar\partial_{B|_\Sigma}^* \\
    \sqrt{2}\partial_{B|_\Sigma} & -i\nabla_{B,\partial_t}
  \end{pmatrix}
  \begin{pmatrix}
    \Psi^+ \\ \Psi^-
  \end{pmatrix} = 0, \label{eq2:SW3-4}
\end{equation}
where $\nabla_{B,\partial_t}$ denotes the $\spinc$ covariant derivative of $B$ evaluated in the $\partial_t$ direction.
Thus, equations (\ref{eq2:SW3-1})--(\ref{eq2:SW3-4}) yield for us the Seiberg-Witten equations on $[0,\infty)\times\Sigma$.

In the same way that the Seiberg-Witten equations on a product $4$-manifold can be written as the downward flow of the Seiberg-Witten vector field induced from the slice $3$-manifold (when the configuration in question is in temporal gauge), we reinterpret the Seiberg-Witten equations on $[0,\infty)\times\Sigma$ as a downward flow of a vector field on the configuration space of $\Sigma$.  To do this, we can consider the oriented $4$-manifold $S^1 \times [0,\infty) \times \Sigma$ and regard configurations on $[0,\infty) \times \Sigma$ as $S^1$ invariant.  If we do this, and we place $(B,\Psi)$ in temporal gauge, then we can regard $(B,\Psi)$ as a downward flow for the Seiberg-Witten vector field on $S^1 \times \Sigma$:
\begin{equation}
  \frac{d}{dt}(B,\Psi) = -SW^{S^1 \times \Sigma}_3\left((B(t),\Psi(t))|_{S^1 \times \Sigma}\right). \label{eq2:CSDflow}
\end{equation}
Here, $SW^{S^1 \times \Sigma}_3$ denotes the gradient of the Chern-Simons-Dirac functional on $S^1 \times \Sigma$.  The Clifford multiplication $\tilde \rho$ on $S^1 \times \Sigma$ is such that\footnote{As in \cite{KM}, given Clifford multiplication $\rho_3$ on a $3$-manifold $Y$, the Clifford multiplication $\rho_4^Y$ on $S^1_\theta \times Y$, regarded as endomorphisms of the self-dual spinor bundle (which is isomorphic to the anti-self-dual spinor bundle in this case) can be chosen so that $\rho_4(\partial_\theta) = \mathrm{id}$ and $\rho_4(v) = \rho_3(v)$ for all $v \in TY$.  This is the natural choice for interpreting the Seiberg-Witten equations on $S^1 \times Y$ as a gradient flow of the Chern-Simons-Dirac functional on $Y$. When $Y = [0,\infty) \times \Sigma$, then on $Y' = S^1_\theta \times \Sigma$, the Clifford multiplication $\rho_3^{Y'}(\cdot) := \rho_4(\partial_t)^{-1}\rho_4(\cdot)$ is the relevant Clifford multiplication to consider since it switches the roles of $\theta$ and $t$, i.e. $\rho_4^{Y'}(\partial_t) = \rho_4^Y(\partial_\theta)$.} $\tilde \rho(\partial_\theta) = \rho(\partial_t)^{-1}$, where $\theta$ is the coordinate on $S^1$, and $\tilde \rho|_{T\Sigma} = \rho(\partial_t)^{-1}\rho|_{T\Sigma}$.  Now for any $S^1$ invariant configuration $(B,\Psi)$ on $S^1 \times \Sigma$, the Chern-Simons-Dirac functional on $S^1 \times \Sigma$ is given by
\begin{equation}
  CSD(B,\Psi) = \frac{1}{2}\int_\Sigma \Re(\Psi,D_{B|_\Sigma}\Psi). \label{eq2:CSD}
\end{equation}
Here, the Chern-Simons term drops out since $B$ has no $S^1$ dependence or $S^1$ component, the operator $D_{B|_\Sigma}$ is the induced Dirac operator on $\Sigma$, and the length of $S^1$ is normalized to unity.\\


\noindent\textbf{Notation. } We write $C$ to denote a connection on $\Sigma$ and $\Upsilon$ to denote a spinor on $\Sigma$, i.e. $(C,\Upsilon)$ is an element of the configuration space $\fC(\Sigma) = \A(\Sigma) \times \Gamma(\S_\Sigma)$ on $\Sigma$.  This is to keep our notation consistent with \cite{KM}, which in general, writes $(B,\Psi)$ for a $3$-dimensional configuration and $(A,\Phi)$ for a $4$-dimensional configuration.  Likewise, we use $c$ to denote a $1$-form on $\Sigma$.  

In light of (\ref{eq2:CSD}), we define the Chern-Simons-Dirac functional $CSD^\Sigma$ on $\fC(\Sigma)$ by \label{p:CSDSigma}
\begin{equation}
  CSD^\Sigma(C,\Upsilon) = \frac{1}{2}\int_\Sigma \Re(\Upsilon, D_C\Upsilon), \qquad (C,\Upsilon) \in \fC(\Sigma),
\end{equation}
where $D_C: \Gamma(\S_\Sigma) \to \Gamma(\S_\Sigma)$ is the $\spinc$ Dirac operator determined from $C$.  The $L^2$-gradient of this functional is given by \label{p:SW2}
\begin{align}
  SW_2(C,\Upsilon) & := \nabla_{(C,\Upsilon)}CSD^\Sigma\\
&= (\tilde \rho_\Sigma^{-1}(\Upsilon\Upsilon^*)_0, D_C\Upsilon),
\end{align}
where $\tilde \rho_\Sigma^{-1}: i\frak{su}(\S_\Sigma) \to T\Sigma$ is the map $\tilde \rho^{-1}: i\frak{su}(\S_\Sigma) \to T(S^1 \times \Sigma)$ composed with the projection onto the $T\Sigma$ factor.  We can consider the formal downward gradient flow of $CSD^\Sigma$ on $\fC(\Sigma)$
\begin{equation}
  \frac{d}{dt}(C,\Upsilon) = -SW_2(C,\Upsilon). \label{eq2:SW2flow}
\end{equation}

Regarding the $S^1$ invariant configuration $(B,\Psi)$ in (\ref{eq2:CSDflow}) as a path of configurations $(B,\Psi) = (C(t),\Psi(t))$ in $\fC(\Sigma)$, we see that (\ref{eq2:SW2flow}) differs from (\ref{eq2:CSDflow}) from the fact that the term $SW_3^{S_1 \times \Sigma}(B,\Psi)$ contains a $d\theta$ component, where $\theta$ denotes the coordinate on $S^1$.  However, because $B$ is $S^1$-invariant and therefore has no $d\theta$ component, equation (\ref{eq2:CSDflow}) implies that the $d\theta$ component of $SW_3^{S_1 \times \Sigma}(B,\Psi)$ is identically zero, i.e., we have a constraint.

Since $\tilde\rho(\partial_\theta) = \rho(\partial_t)^{-1}$, this constraint is none other than the equation (\ref{eq2:SW3-1}).
In light of this, given $(C,\Upsilon) \in \fC(\Sigma)$, define the map \label{p:momentmap}
\begin{equation}
  \mu(C,\Upsilon) = \check{*}F_C + \frac{i}{2}(|\Upsilon_-|^2 - |\Upsilon_+|^2). \label{eq2:mu}
\end{equation}
Here $\check{*}$ is the Hodge star on $\Sigma$ and $\Upsilon = (\Upsilon_+,\Upsilon_-)$ is the decomposition of $\Upsilon \in \Gamma(\S_\Sigma)$ induced by the splitting (\ref{eq2:Siso}).

Recall that the gauge group $\G(\Sigma) = \Maps(\Sigma,S^1)$ acts on $\fC(\Sigma)$ via
$$(C,\Upsilon) \mapsto g^*(C,\Upsilon) = (C - g^{-1}dg, g\Upsilon), \qquad g \in \G(\Sigma).$$
We have the following proposition concerning the map $\mu$:

\begin{Proposition}\label{PropMomentMap} $\,$
\begin{enumerate}
  \item The map $\mu: \fC(\Sigma) \to \Omega^0(\Sigma; i\R)$ is the moment map for $\fC(\Sigma)$ associated to the gauge group action of $\G(\Sigma)$.  Here, the symplectic form on $\fC(\Sigma)$ is given by
      $$\omega((a,\phi),(b,\psi)) = \int_\Sigma a \wedge b + \int_\Sigma \Re(\phi,\rho(-\partial_t)\psi), \qquad (a,\phi), (b,\psi) \in \T_\Sigma.$$
  \item If $\Upsilon \not\equiv 0$, then $d_{(C,\Upsilon)}\mu: T_{(C,\Upsilon)}\fC(\Sigma) \to \Omega^0(\Sigma;i\R)$ is surjective.
  \item A configuration $(B,\Psi) = (C(t),\Upsilon(t))$ in temporal gauge on $[0,\infty)\times\Sigma$ solves $SW_3(B,\Psi) = 0$ if and only if $(C(t),\Upsilon(t))$ solves
      \begin{align}
        \frac{d}{dt}(C(t),\Upsilon(t)) &= -SW_2(C(t),\Upsilon(t)) \label{eq2:SW2flow1}\\
        \mu(C(t),\Upsilon(t)) &= 0, \qquad t > 0. \label{eq2:SW2flow2}
      \end{align}
  \item Suppose we consider the perturbed equations $SW_3(B,\Psi) = \eta$, where $\eta = \eta^0 dt$ and $\eta^0 \in \Omega^0(\Sigma; i\R)$ is time-independent.  Then (iii) holds but with (\ref{eq2:SW2flow2}) replaced by
      \begin{equation}
        \mu(C(t),\Upsilon(t)) = \eta^0, \qquad t > 0. \label{eq2:SW2flow2'}
      \end{equation}
   \item For any $(C,\Upsilon) \in \fC(\Sigma)$, we have $d_{(C,\Upsilon)}\mu\Big(SW_2(C,\Upsilon)\Big) = 0$, that is, $SW_2\cc$ is tangent to any level set of $\mu$.
\end{enumerate}
\end{Proposition}

\Proof (i) This is the statement that at every $\cc \in \fC(\Sigma)$, every $(c,\upsilon) \in T_\cc\fC(\Sigma)$, and every $\xi \in \Omega^0(\Sigma; i\R)$, we have
$$\int_\Sigma d_\cc\mu(c,\upsilon) \cdot \xi = \omega((-d\xi,\xi\Upsilon),(c,\upsilon)).$$
Verifying this is a straightforward computation.

(ii) The range of $\check *d: \Omega^1(\Sigma; i\R) \to \Omega^0(\Sigma; i\R)$ consists of precisely those functions that integrate to zero on $\Sigma$.  Suppose $f$ is orthogonal to the image of $d_\cc\mu$.  If $\Upsilon \not\equiv 0$, then one can find $\upsilon \in \Gamma(\S_\Sigma)$ such that $d_\cc\mu(0,\upsilon) = \int_\Sigma f$.  It then follows that one can find a $1$-form $c$ such that $d_\cc(c,\upsilon) = f$.

Statements (iii) and (iv) follow from (\ref{eq2:SW3-1}).

(v) First observe that since $\mu$ is the moment map for the gauge group action on $\fC(\Sigma)$, the kernel of its differential is the symplectic annihilator of the tangent space to the gauge orbit:
\begin{align*}
  \ker d_\cc\mu &= \{(-d\xi,\xi\Upsilon) : \xi \in \Omega^0(\Sigma; i\R)\}^{\bot_\omega}\\
  &= J\{(-d\xi,\xi\Upsilon) : \xi \in \Omega^0(\Sigma; i\R)\}^{\bot}.
\end{align*}
Here,
$$J := (-\check{*}, \rho(\partial_t)): \Omega^1(\Sigma; i\R) \oplus \Gamma(\S_\Sigma) \circlearrowleft$$
is the compatible complex structure for $\omega$ (where the associated inner product is the usual $L^2$ inner product).  Thus, to show $d_\cc\mu(SW_2\cc) = 0$, it suffices to show that $J\cdot SW_2\cc$ is perpendicular to the tangent space to the gauge orbit of $\cc$.  For this, it suffices to show that $J\cdot SW_2\cc$, like $SW_2\cc$, is the gradient of a gauge-invariant functional.  A simple computation shows that the gradient of the functional
\begin{equation}
  \cc \mapsto \frac{1}{2}\int_\Sigma(\Upsilon,\rho(\partial_t)D_C\Upsilon) \label{Ham}
\end{equation}
is $JSW_2\cc$.  Here, we use the fact that, by convention of our choice of Clifford multiplication, $\rho(-\partial_t)\rho(d\Sigma) = 1$ and so $\rho(\partial_t)\rho(c) = \rho(d\Sigma)\rho(c) = \rho(\check{*}c).$\End

The last statement of the above lemma implies that the restriction of the gradient of $CSD^\Sigma$ to the level set $\mu^{-1}(\eta^0)$ is equal to the gradient of $CSD^\Sigma|_{\mu^{-1}(\eta^0)}$ (at points where $\mu^{-1}(\eta^0)$ is smooth).  If $ c_1(\s) \neq \frac{i}{\pi}[\check{*}\eta^0]$, then by (iii), $\mu^{-1}(0)$ is always a smooth submanifold of $\fC(\Sigma)$ since it contains no flat connections.  Thus, we have the following corollary:

\begin{Corollary}\label{CorSW2Flow}
Suppose $c_1(\s) \neq \frac{i}{\pi}[\check{*}\eta^0]$.  Then $\mu^{-1}(\eta^0)$ is a smooth submanifold of $\fC(\Sigma)$ and modulo gauge, solutions to $SW_3(B,\Psi) = 0$ on $Y$ correspond to (formal) downward gradient flow lines of $CSD^\Sigma|_{\mu^{-1}(\eta^0)}$.
\end{Corollary}

\noindent Thus, the bulk of our analysis consists in understanding the gradient flow of $CSD^\Sigma|_{\mu^{-1}(\eta^0)}$. As a remark, note that our flow is also Hamiltonian with respect to the functional (\ref{Ham}).

\subsection{The Vortex Equations}

Let $d = \frac{1}{2}\left<c_1(\s),[\Sigma]\right>$.  From now on, we always assume $\eta^0$ is chosen such that
$$\frac{i}{2\pi}\int_\Sigma \eta^0 \neq d,$$
so that $\mu^{-1}(\eta^0)$ is a smooth manifold.

Our first task is to understand the set of critical points of $CSD^\Sigma|_{\mu^{-1}(\eta^0)}$.  We have the following two facts. First, the critical points of $CSD^\Sigma|_{\mu^{-1}(\eta^0)}$ have an explicit description in terms of the space of vortices on $\Sigma$.  Second, this critical set is Morse-Bott nondegenerate with respect to $CSD^\Sigma|_{\mu^{-1}(\eta^0)}$.  This is in contrast to the case $\frac{i}{2\pi}\int_\Sigma \eta^0 = d$, where although the critical set of $CSD^\Sigma|_{\mu^{-1}(0)}$ is just the space of flat connections on $\Sigma$, this set is in general Morse-Bott degenerate.\footnote{A flat connection $C$ will be Morse-Bott degenerate precisely when $\ker D_C \neq 0$.}

For the sake of completeness, we describe in detail the correspondence between the critical set of $CSD^\Sigma|_{\mu^{-1}(\eta^0)}$ and the space of vortices, following \cite{MOY}.  Recall that the vortex equations on $\Sigma$ are given by the following.  Given a line bundle $E \to \Sigma$ over $\Sigma$ of degree $k$, a Hermitian connection $A$ on $E$, a section $\psi \in \Gamma(E)$, and a function $\tau \in \Omega^0(\Sigma; i\R)$, the vortex equations are given by
\begin{align}
  \check{*}F_A - \frac{i|\psi|^2}{2} &= \tau \label{eq2:vort1}\\
  \bar\partial_A\psi & = 0. \label{eq2:vort2}
\end{align}
Here, $\bar\partial_A: E \to K_\Sigma^{-1}\otimes E$ is the holomorphic
structure on $E$ determined by $A$.  Observe that if $k > \frac{i\tau}{2\pi}$, there are no solutions to (\ref{eq2:vort1})--(\ref{eq2:vort2}) by a simple application of the Chern-Weil theorem.  When $0 < k < \int_\Sigma \frac{i\tau}{2\pi}$, then by \cite{GP}, the moduli space of gauge equivalence classes of solutions  \label{p:fV} $\V_{k,\tau}(\Sigma)$ to (\ref{eq2:vort1})--(\ref{eq2:vort2}) can be naturally identified with the space of effective divisors of degree $k$ on $\Sigma$, i.e. the $k$-fold symmetric product $\Sym^k(\Sigma)$ of $\Sigma$.  This identification is given by mapping a solution $(A,\psi)$ to the set of zeros of the (nontrivial) holomorphic section $\psi$. Because $\Sym^k(\Sigma)$ is independent of $\tau$, we will often simply denote the moduli space of degree $k$ vortices by $\V_k(\Sigma)$.  Likewise, we will denote the space of solutions (without dividing by gauge) to (\ref{eq2:vort1})--(\ref{eq2:vort2}) by $\fV_k(\Sigma)$.

Observe that if $k < 0$, one may instead consider the equations
\begin{align}
  \check{*}F_A + \frac{i|\psi|^2}{2} &= -\tau \label{eq2:antivort1}\\
  \partial_A\psi & = 0, \label{eq2:antivort2}
\end{align}
which become equivalent to (\ref{eq2:vort1})--(\ref{eq2:vort2}) via complex conjugation.  We will call the equations (\ref{eq2:antivort1})--(\ref{eq2:antivort2}) the anti-vortex equations.  Thus, the  space of solutions to (\ref{eq2:antivort1})--(\ref{eq2:antivort2}), which we denote by  $\bar\fV_{k,\tau}(\Sigma)$, is nonempty for $\int_\Sigma -\frac{i\tau}{2\pi} < k < 0$ and, its moduli space of gauge equivalence classes $\bar\V_{k,\tau}(\Sigma)$, can can be identified with $\V_{|k|,\tau}(\Sigma)$.

The equations that determine the critical points of $CSD^\Sigma$ which belong to the zero set of the moment map are given by $\mu(C,\Upsilon) = 0$ and $SW_2(C,\Upsilon) = 0$.  More explicitly, these equations are given by
\begin{align}
  \check{*}F_C + \frac{i}{2}(|\Upsilon_-|^2-|\Upsilon_+|^2) &= \eta^0 \label{eq2:SW2crit1} \\
  \Upsilon_+\bar\Upsilon_- & = 0 \label{eq2:SW2crit2} \\
  \bar\Upsilon_+\Upsilon_- & = 0 \label{eq2:SW2crit3} \\
  \bar\partial_C\Upsilon_+ &= 0 \label{eq2:SW2crit4} \\
  \bar\partial_C^*\Upsilon_- &= 0. \label{eq2:SW2crit5}
\end{align}

We can now see the correspondence between equations (\ref{eq2:SW2crit1})--(\ref{eq2:SW2crit5}) and the vortex equations (\ref{eq2:vort1})--(\ref{eq2:vort2}).  Equations (\ref{eq2:SW2crit4})--(\ref{eq2:SW2crit5}) and unique continuation for  Dirac operators imply that (\ref{eq2:SW2crit2})--(\ref{eq2:SW2crit3}) forces
$$\Upsilon_+ \equiv 0 \qquad \mathrm{or} \qquad \Upsilon_- \equiv 0.$$

Let $g$ denote the genus of $\Sigma$. (If $\Sigma$ is not connected, we can of course work on each component of $\Sigma$ separately.) Pick a connection $C_\Sigma$ on $K_\Sigma^{1/2}$ and define $\tau = \check{*}F_{C_\Sigma}$.  
Let $\bar C_\Sigma$ denote the corresponding dual connection on $K_\Sigma^{-1/2}$.  Then (\ref{eq2:SW2crit1}) is equivalent to each of the following equations
\begin{align}
  \check{*}F_{C \otimes C_\Sigma} + \frac{i}{2}(|\Upsilon_-|^2-|\Upsilon_+|^2) & = \tau + \eta^0 \label{eq2:SWvort1} \\
  \check{*}F_{C \otimes \bar C_\Sigma} + \frac{i}{2}(|\Upsilon_-|^2-|\Upsilon_+|^2) & = -\tau + \eta^0 \label{eq2:SWvort2}
\end{align}
Given $\Sigma$ and $d$, let
\begin{equation}
  k_\pm = \pm k_g + d, \label{kpm}
\end{equation}
where $k_g = \deg(K_\Sigma^{1/2})$.  Using the constraints on $k$ for when the vortex and anti-vortex moduli spaces $\V_{k}(\Sigma)$ and $\V_{-|k|}(\Sigma)$ are nonempty, and the fact that a line bundle can have nontrivial holmorphic sections only if it has nonnegative degree, it is easy to see that the following situation holds:

\begin{Lemma}\label{LemmaVortices}
With notation as above, we have the following:
\begin{enumerate}
  \item Suppose $-k_g \leq d < \frac{i}{2\pi}\int_\Sigma \eta^0$.  Then the space of critical points of $CSD^\Sigma|_{\mu^{-1}(\eta^0)}$ corresponds precisely to the space of vortices $\fV_{k_+}(\Sigma)$ under the correspondence $(C,\Upsilon) \mapsto (C \otimes C_\Sigma, \Upsilon_+)$.  Here $\Upsilon_-$ vanishes identically.
  \item Suppose $\frac{i}{2\pi}\int_\Sigma \eta^0 < d \leq k_g$.  Then the space of critical points of $CSD^\Sigma|_{\mu^{-1}(\eta^0)}$ corresponds precisely to the space of anti-vortices $\bar\fV_{k_-}(\Sigma) \cong \fV_{|k_-|}(\Sigma)$ under the correspondence $(C,\Upsilon) \mapsto (C \otimes \bar C_\Sigma, \Upsilon_-)$. Here, $\Upsilon_+$ vanishes identically.
  \item For $\frac{i}{2\pi}\int_\Sigma \eta^0 = d$, the space of critical points of $CSD^\Sigma|_{\mu^{-1}(\eta^0)}$ corresponds precisely to the space of flat connections on $\Sigma$ under the correspondence $C \mapsto C \otimes \bar C_\Sigma$.  Here, $\Upsilon$ vanishes identically.
  \item For all other values of $d$, the set of critical points of $CSD^\Sigma|_{\mu^{-1}(\eta^0)}$ is empty.
  \item For all choices of $\eta^0$ as above, except in case (iii), the critical set of $CSD^\Sigma|_{\mu^{-1}(\eta^0)}$ is Morse-Bott nondegenerate.
\end{enumerate}
\end{Lemma}

\Proof Statements (i)-(iv) follow from the preceding analysis.  

We need only prove (v). This amounts to showing the following.  Given any configuration $\cco \in \fC(\Sigma)$, let
\begin{equation}
  \H_{2,\cco}: \T_\Sigma \to \T_\Sigma, \label{eq2:Hess2def}
\end{equation}
denote the \textit{Hessian} of $CSD^\Sigma$ at $\cco$, which is the operator obtained by linearizing the map $SW_2: \fC(\Sigma) \to \T_\Sigma$ at $\cco$.  If $\cco$ is a vortex, we need to show that the restricted operator \label{p:Hess2}
\begin{equation}
  \H_{2,\cco}: T_\cco\mu^{-1}(0) \to T_\cco\mu^{-1}(0), \label{eq2:Hess2}
\end{equation}
has kernel equal to precisely the tangent space to the space of vortices at $\cco$.  Without loss of generality, suppose we are in case (i).  Then if we linearize the equations (\ref{eq2:SW2crit1})--(\ref{eq2:SW2crit5}) at a vortex, then since $\Upsilon_- \equiv 0$ and $\Upsilon_+$ vanishes only on a finite set of points, unique continuation shows that an element of the kernel of the linearized equations must have vanishing $\Psi_-$ component.  It follows that the only nontrivial equations we obtain are those obtained from linearizing (\ref{eq2:SW2crit1}) and (\ref{eq2:SW2crit4}), which yields for us precisely the linearization of the vortex equations. On the other hand, the space of vortices are cut out transversally by the vortex equations. (This is because the set $\{(A,\psi): \bar\partial_A\psi = 0, \psi \not\equiv 0\}$ is a gauge-invariant Kahler submanifold of $\fC(\Sigma)$, the left-hand side of (\ref{eq2:vort1}) is the moment map for this submanifold, and the gauge group acts freely on this submanifold.) It follows that the kernel of the map $\H_{2,\cco}$ above is precisely the tangent space to the space of vortices.  This finishes the proof of Morse-Bott nondegeneracy.\End

Given the above lemma, we abuse notation by letting $\fV_k(\Sigma)$ denote the set of critical points of $CSD^\Sigma|_{\mu^{-1}(\eta^0)}$, with $k$ suitably defined as given by the above.   \textit{We always assume that we are in the Morse-Bott situation hereafter.}  We write $\V_k(\Sigma)$ to denote the quotient of $\fV_k(\Sigma)$ by the gauge group, and it can be identified with $Sym^k(\Sigma)$.  Observe that for $\eta^0$ such that $\frac{i}{2\pi}\int_\Sigma \eta^0$ is not an integer, $d$ and hence $k$ can be arbitrary integers (the moduli spaces being empty for negative $k$), so that our analysis works for arbitrary $\spinc$ structures.

For every $k$, note that the symplectic form on $\fC(\Sigma)/\!/\G(\Sigma)$ restricts to a symplectic form on the vortex moduli space $\V_k(\Sigma)$.  We will refer to elements of either $\fV_k(\Sigma)$ or $\V_k(\Sigma)$, for any $k$, simply as vortices.  When $\Sigma$ and $k$ are fixed, we will often write $\fV$ and $\V$ for brevity.

\subsection{The Flow on a Slice}

We now assume $k$ is fixed.  In order to place ourselves in an elliptic situation and in a situation where we can apply Morse-Bott estimates to our configurations, we have to choose the right gauge for our equations.  As it turns out, choosing a suitable gauge requires some careful setup.  Our work here is modeled off that of \cite{MMR}, which studies the flow one obtains for the instanton equations on a cylindrical $4$-manifold. To describe the gauge fixing procedure, we recall the basic gauge theoretic decompositions of the configuration space on $\Sigma$ and its tangent spaces.

Our analysis proceeds similarly to the case of a closed $3$-manifold \cite{KM}.  Given a configuration $\cc \in \fC(\Sigma)$, define
\begin{align*}
  \T_{(C,\Upsilon)} = T_{(C,\Upsilon)}\fC(\Sigma) = \Omega^1(\Sigma;i\R) \oplus\Gamma(\S_\Sigma)
\end{align*}
to be the tangent space to $(C,\Upsilon)$ of $\fC(\Sigma)$. If the basepoint is unimportant, we write $\T_\Sigma$ for any such tangent space.  The infinitesimal action of the gauge group on $\fC(\Sigma)$ leads us to consider the following operators
\begin{align*}
  \bd_\cc: \Omega^0(\Sigma; i\R) & \to \T_\Sigma\\
  \xi & \mapsto (-d\xi,\xi\Upsilon)\\
  \bd_\cc^*: \T_\Sigma & \to \Omega^0(\Sigma; i\R)\\
  (c,\upsilon) & \mapsto -d^*c + i\Re(i\Upsilon,\upsilon).
\end{align*}
From these operators, we obtain the following subspaces of $\T_{(C,\Upsilon)}$, which are the tangent space to the the gauge orbit through $(C,\Upsilon)$ and its orthogonal complement, respectively:
\begin{align*}
  \J_{(C,\Upsilon)} &= \im \bd_\cc\\
  \K_{(C,\Upsilon)} &= \ker \bd_\cc^*.
\end{align*}

As usual, we must consider the Banach space completion of the configuration spaces and the above vector spaces. Given a manifold $M$ and $s \geq 0$, let $H^s(M)$ denote the Sobolev space of functions that have $s$ (fractional) derivatives belonging to $L^2(M)$. Write $\fC^s(\Sigma)$ to denote the $H^s(\Sigma)$ completion of the configuration space on $\Sigma$.  Its tangent spaces are  isomorphic to $\T^s_\Sigma$, the $H^s(\Sigma)$ completion of $\T_\Sigma$. For sufficiently regular $\cc$, we obtain the following subspaces of $\T^s_\Sigma$:
\begin{align*}
  \J_{(C,\Upsilon)}^s &= \{\bd_\cc\xi  : \xi \in H^{s+1}\Omega^0(\Sigma;i\R)\}\\
  \K_{(C,\Upsilon)}^s &= \{(c,\upsilon) \in \T_\cc^s: \bd_\cc^*(c,\upsilon) = 0\}.
\end{align*}

We have the following gauge-theoretic decompositions of the tangent space and configuration space:

\begin{Lemma} \label{LemmaSliceConfig}
Let $s > 0$.
  (i) Then for any $\cc \in \fC^s(\Sigma)$, we have an $L^2$ orthogonal decomposition
  \begin{equation}
    \T_\cc^s = \J_\cc^s \oplus \K_\cc^s. \label{eq2:T=JK}
  \end{equation}
  (ii) Define the slice
  $$\fS_\cco^s := \cco + \K_\cco^s$$
  through $\cco$ in $\T_\cc^s$. There exists an $\eps > 0$ such that if $\cc \in \fC^s(\Sigma)$ satisfies $\|\cc - \cco\|_{H^s(\Sigma)} < \eps$, then there exists a gauge transformation $g \in \G^{s+1}(\Sigma)$ such that $g^*\cc \in \fS_\cco^s$ and $\|g^*\cc - \cco\|_{H^s(\Sigma)} \leq c_s\|\cc - \cco\|_{H^s(\Sigma)}$.
\end{Lemma}

\Proof (i) This lemma follows from standard elliptic theory, cf. \cite[Lemma 3.4]{N1}.  (ii) This is an immediate consequence of the inverse function theorem and the fact that $\fS_\cco^s$ is a local slice for the gauge action.\End

For $s > 0$, define the quotient configuration space \label{p:fB}
$$\fB^s(\Sigma) = \fC^s(\Sigma)/\G^{s+1}(\Sigma).$$
Away from the reducible configurations (i.e. those for which the spinor vanishes identically), this quotient space is a Hilbert manifold modeled on the above local slices (see \cite{KM}). The decomposition (\ref{eq2:T=JK}) allows us to define the complementary projections $\Pi_{\J_\cc^s}$ and $\Pi_{\K_\cc^s}$ of $\T_\cc^s$ onto $\J_\cc^s$ and $\K_\cc^s$, respectively.

Let us return to the smooth setting for the time being.  Denote the quotient of the smooth configuration space by the smooth gauge group by
$$\fB(\Sigma) = \fC(\Sigma)/\G(\Sigma).$$
Our first task is to rewrite the Seiberg-Witten equations on $Y$ in a suitable gauge when the monopole in question is close to a vortex.  This is so that we may exploit the Morse-Bott nature of the critical set, which we perform in the next section.  \\

\noindent\textbf{Notation. } \label{p:gamma} To simplify notation a bit, and to make it bear similarity with that of the standard reference \cite{KM}, we introduce the following notation.  We will write $\frak{a}$ to denote a critical point of $CSD^\Sigma|_{\mu^{-1}(\eta^0)}$, i.e. a vortex. We will always assume $\fa$ is smooth, unless otherwise stated, since this can always be achieved via a gauge transformation. Given a configuration $(B,\Psi)$ on $Y = [0,\infty) \times \Sigma$, we can write it as
$$(B,\Psi) = (C(t) + \beta(t)dt, \Upsilon(t))$$
where $(C(t),\Upsilon(t))$ is a path of configurations in $\fC(\Sigma)$ and $\beta(t)$ is a path in $\Omega^0(\Sigma;i\R)$.  As shorthand, we will often write $\g$ for the configuration $(B,\Psi)$ and $\check \g(t)$ for the path $(C(t),\Upsilon(t))$.  Given a vortex $\fa$, we write $\g_\fa \in \fC(Y)$ to denote the time-translation invariant path identically equal to $\fa$.\\

To simplify notation in what follows, we temporarily drop the superscript $s$ and work with smooth objects since everything works mutatis mutandis in the $H^s$ topology. Given any vortex $\fa \in \fV$, define
$$CSD^\Sigma_\fa = CSD^\Sigma|_{\fS_\fa} - CSD^\Sigma(\fa).$$
to be the restriction of $CSD^\Sigma$ to the slice $\fS_\fa$ of smooth configurations through $\fa$, normalized by a constant for convenience.  Note that $CSD^\Sigma$ has a constant value on its critical set, since it is connected.

Since $SW_2\cc = \nabla_\cc CSD^\Sigma$ is the gradient of the gauge-invariant functional $CSD^\Sigma$, we know that $\nabla_\cc CSD^\Sigma$ is orthogonal to $\J_\cc$ and hence lies in $\K_\cc$.  On the other hand, if $\cc \in \fS_\fa$, then the gradient of $CSD^\Sigma_\fa$ satisfies
$$\nabla_{\cc}CSD^\Sigma_\fa \in \K_{\fa},$$
since a priori, this gradient must be tangent to the slice.  For $\cc$ close enough to $\a$, then the space $J_\cc$, which is automatically complementary to $\K_\cc$, is also complementary to $\K_\a$, and so $\nabla_\cc CSD^\Sigma$ and $\nabla_{\cc}CSD^\Sigma_\fa$ differ by an element of $\J_\cc$.  This suggests we introduce the following inner product structure on the tangent bundle of a neighborhood $\fS_\a(\delta)$ of the slice (instead of the usual $L^2$ inner product).  Namely, mimicking the construction in \cite{MMR}, consider the inner product
\begin{equation}
  \left<x,y\right>_{\a,\cc} := \left(\Pi_{\K_\cc}x, \Pi_{\K_\cc}y\right)_{L^2(\Sigma)}, \qquad x,y \in T_\cc\fS_\a(\delta) \label{eq2:innerprod}
\end{equation}
where $\left(\cdot,\cdot\right)_{L^2(\Sigma)}$ is the usual $L^2$ inner product on $\T_\Sigma$, and $\Pi_{\K_\cc}$ is the orthogonal projection onto $\K_\cc$ with kernel $\J_\cc$.  As noted, for $\cc$ sufficiently close to $\a$, the map $\Pi_{\K_\cc}: \K_\a \to \K_\cc$ is an isomorphism.  Specifically, by the same analysis as in \cite[Remark 4.3]{N1}, $\cc$ in a small $H^{1/2}(\Sigma)$ ball $U$ around $\a$ is sufficient.  Observe that the inner product $\left<\cdot,\cdot\right>_{\a,\cc}$ naturally arises from pulling back the $L^2$ inner product on the irreducible part of the quotient configuration space $\fC(\Sigma)/\G(\Sigma)$.

Then if we endow the neighborhood $U$ with the inner product (\ref{eq2:innerprod}), we can explicitly write $\nabla_{\cc}CSD^\Sigma_\fa$ as follows.  Let $\Pi_{\K_\a,\J_\cc}$ denote the projection onto $\K_\a$ through $\J_\cc$, which exists for $\cc \in U$ and $U$ sufficiently small.  Then
\begin{equation}
  \nabla_{\cc}CSD^\Sigma_\fa = \Pi_{\K_\a,\J_\cc}SW_2\cc, \label{eq2:PiKJSW2}
\end{equation}
or in other words, there exists a well-defined map
\begin{equation}
  \Theta_\a: U \to \Omega^0(\Sigma; i\R)
\end{equation}
such that
\begin{equation}
  \nabla_{\cc}CSD^\Sigma_\fa = SW_2\cc - \bd_\cc\Theta_\a\cc. \label{eq2:thetadecomp}
\end{equation}
(The map $\Theta_\a$ is well-defined since the operator $\bd_\cc$ is injective for $\cc$ irreducible, which holds for $U$ small.)  The decomposition (\ref{eq2:thetadecomp}) is important because it relates the gradient vector field $\nabla_{\cc}CSD^\Sigma_\fa$ to the vector field $SW_2\cc$ by an infinitesimal action of the gauge group at the configuration $\cc$.  (Had we used the usual $L^2$ inner product, the analogous ansatz would have yielded an infinitesimal action of the gauge group at $\a$ instead of the configuration $\cc$ in question.)

Borrowing the terminology of \cite{MMR}, we introduce the following definition:

\begin{Definition}Fix $s \geq 1/2$ and $\eta^0 \in \Omega^0(\Sigma; i\R)$.
\begin{enumerate}
  \item For any smooth vortex $\fa$, any open subset of $\fS_\fa^s \cap \mu^{-1}(\eta^0)$ of the form
  $$U_\fa(\delta) := \{\cc \in \fS_\fa^s \cap \mu^{-1}(\eta^0) : \|\cc - \a\|_{H^{1/2}} < \delta\}$$
   for some small $\delta > 0$ is said to be a \textit{coordinate patch} at $\fa$.  We will often write $U_\fa$ to denote any such coordinate patch.  We always assume that the (sufficiently small) coordinate patch $U_\fa$ is endowed with the inner product (\ref{eq2:innerprod}) on its tangent bundle.
  \item Let $I$ be a subinterval of $[0,\infty)$. Given a coordinate patch $U_\fa$ about a vortex $\fa$, we say that a configuration $\g \in \fC([0,\infty) \times \Sigma)$ is in \textit{standard form} on $I \times \Sigma$ with respect to $U_\fa$ if $\cg(t) \in U_\fa$ for all $t \in I$.
\end{enumerate}
\end{Definition}
Our choice of defining $H^{1/2}$ open neighborhoods comes from our energy analysis of the next section. The value of $s$ is immaterial for now and can be assumed as large as desired ($s \geq 2$ is sufficient).  

The upshot of the above formalism is the following.  Given a path of configurations $\cct$ that is sufficiently near a vortex $\a$ for all time $t$, we can gauge fix this path so that the new path lies in some neighborhood of $\a$ in the slice $\fS_\a$ for all time.  The relevant situation is when this path of configurations is a monopole on $Y = [0,\infty) \times \Sigma$ in temporal gauge.  When we perform such a gauge-fixing, two things happen. First, the resulting configuration $\g$ determines a path $\cg(t)$ in a coordinate patch $U_\a$ (i.e., it is in standard form), since our perturbed monopole always determines a path in the appropriate level set of the moment map by Corollary \ref{CorSW2Flow}.  Second, $\g$ is no longer in temporal gauge. Nevertheless, the next lemma tells us that the resulting configuration $\g$ is completely determined by the path $\cg(t)$.  Moreover, the path $\cg(t)$ is simply a gradient flow line of $CSD^\Sigma_\a$ restricted to $U_\a$.

\begin{Lemma} \label{LemmaStandFlow}
Let $\fa$ be a vortex and $U_\fa$ a coordinate patch.  Let $\g = (C(t) + \beta(t)dt, \Upsilon(t))$ be a configuration on $[0,\infty) \times \Sigma$ in standard form on $[T_0,T_1] \times \Sigma$ with respect to $U_\a$.  Then $\g$ satisfies $SW_3(\g) = (\eta^0dt, 0)$ if and only if
\begin{equation}
\begin{split}
  \frac{d}{dt}\check\g(t) & = -\nabla_{\cg(t)} CSD_\a  \\
  \bd^*_\a(\cg(t) - \fa) &= 0\\
  \mu(\cg(t)) &= \eta^0,\\
  \beta(t) &= \Theta_\fa(\check\g(t)), \qquad T_0 < t < T_1.
\end{split}\label{eq2:stdform}
\end{equation}
\end{Lemma}

\Proof The monopole equations $SW_3(\g) = (\eta^0dt, 0)$, as given by (\ref{eq2:SW3-4}), are precisely
\begin{equation}
\begin{split}
  \frac{d}{dt}\check\g(t) & = -SW_2(\cg(t)) + \bd_{\check\g(t)}\beta(t)  \\
  \mu(\cg(t)) &= \eta^0.
\end{split} \label{eq2:SW3arb}
\end{equation}
Thus, any solution to (\ref{eq2:stdform}) yields a solution to (\ref{eq2:SW3arb}).  Conversely, suppose we have a solution $\g$ to (\ref{eq2:SW3arb}).  Since $\g$ is in standard form, it satisfies the second equation of (\ref{eq2:stdform}), and taking a time-derivative of this equation, we obtain
$$\bd_\a^*\frac{d}{dt}\cg(t) = 0.$$
The first equation now implies
$$-\bd_\a^*SW_2(\cg(t)) + \bd_\a^*(\bd_{\cg(t)}\beta(t)) = 0.$$
From the definitions, this implies $\beta(t) = \Theta_\a(\cg(t))$.  We now see that $\g$ solves (\ref{eq2:stdform}).\End

We now use this lemma to study the asymptotic behavior of monopoles at infinity.


\section{Asymptotic Convergence and Exponential Decay}

Lemma \ref{LemmaStandFlow} tells us that a solution to the Seiberg-Witten equations on $[0,\infty)\times\Sigma$ in standard form with respect to a small coordinate patch $U_\a$ of a vortex $\a$ satisfies the system of equations (\ref{eq2:stdform}).  These equations tell us that the solution $\g$ is determined by the evolution of the path $\cg(t)$ in $U_\a$, since the normal component $\beta(t)dt$ is determined from $\cg(t)$.  The path $\cg(t)$ is a downward gradient flow for the functional $CSD_\a^\Sigma$ on the coordinate chart $U_\a$, where this latter space has been endowed with the inner product (\ref{eq2:innerprod}).  It is on a sufficiently small coordinate patch $U_\a$ that we can apply standard Morse-Bott type estimates for the function $CSD_\a^\Sigma$.  These estimates imply that any trajectory $\cg(t)$ that stays within $U_\a$ for all time must converge exponentially fast to a critical point.  Moreover, we can deduce that the $L^2(\Sigma)$ length of the path $\cg(t)$ is bounded by the energy of the path, see (\ref{eq2:length}). Here, the energy of a monopole $\g$ is the quantity \label{p:energy}
\begin{equation}
  \E(\g) = \int_0^\infty\|SW_2(\cg(t))\|_{L^2(\Sigma)}^2dt. \label{eq:energy}
\end{equation}
Likewise we can define the energy $\E_I(\g)$ of a configuration on $I \times \Sigma$, for any interval $I = [t_1,t_2]$.  On any such interval for which the energy is finite, the energy is equal to the drop in the Chern-Simons-Dirac functional on $\Sigma$:
$$CSD^\Sigma(\cg(t_1)) - CSD^\Sigma(\cg(t_2)) = \int_{t_1}^{t_2}\|SW_2(\cg(t))\|^2_{L^2(\Sigma)}dt.$$
This is a simple consequence of the fact that a monopole on $I \times \Sigma$ is simply a downward gradient flow line of $CSD^\Sigma$.

Regarding a monopole $(B,\Psi)$ on $I \times \Sigma$ as an $S^1$ invariant configuration on $S^1 \times I \times \Sigma$, with $I$ a compact interval, then we have the following energy identity (see \cite{KM}):
$$CSD^\Sigma(\cg(t_1)) - CSD^\Sigma(\cg(t_2)) = \int_{I \times \Sigma} \left(\frac{1}{4}|F_{B}|^2 + |\nabla_B\Psi|^2 + \frac{1}{4}\Big(|\Psi|^2 + (s/2)\Big)^2 - \frac{s^2}{16}\right)$$
where $s$ is the scalar curvature of $I \times \Sigma$. Thus, modulo gauge, the energy of a monopole controls its $H^1$ norm on finite cylinders.

A key step in understanding the moduli space of finite energy monopoles is to show that if a monopole $\g$ has small enough energy, then there is a vortex $\fa$ and a gauge transformation $g$ on $[0,\infty)\times\Sigma$ such that $g^*\g$ determines a path that stays within some coordinate patch of $\fa$ for all time.  In this way, one can see at an intuitive level what the moduli space of monopoles on $[0,\infty)\times\Sigma$ with small finite energy is.  It is simply a neighborhood of the stable manifold to the space of vortices in the symplectic reduction $\mu^{-1}(\eta^0)/\G(\Sigma)$.  There is some analytic care that must be taken to establish this picture, however, since the coordinate patches we consider only contain $H^{1/2}(\Sigma)$ neighborhoods of a vortex $\a$, whereas the important length estimate (\ref{eq2:length}) is only an $L^2(\Sigma)$ bound.  Nevertheless, it turns out that one can bootstrap the $L^2(\Sigma)$ convergence of the configuration to show that it converges in $H^s(\Sigma)$ exponentially fast to a vortex within a fixed coordinate chart, for all $s \geq 0$.

We begin with the following fundamental estimates for configurations with small energy.  Given any $I$, we write $\fV_I \subset \fC(I \times \Sigma)$ for the space of time translation invariant elements on $I \times \Sigma$ that belong to the space of vortices $\fV = \V(\Sigma)$ for all time.

\begin{Lemma}\label{LemmaSmallEnergy}
We have the following:
\begin{enumerate}
  \item Given a bounded interval $I$, for every gauge invariant neighborhood $V$ of $\fV_I$ in $\fC^1(I \times \Sigma)$, there exists an $\eps > 0$ such that if $\g$ is any monopole on $I\times\Sigma$ satisfying the small energy condition $\int_I \|SW_2(\cg(t))\|_{L^2(\Sigma)}^2 < \eps$, then there exists a gauge transformation $g$ such that $g^*\g \in V$.
  \item For every gauge invariant neighborhood $V_\Sigma$ of $\fV$ in $\fC^1(\Sigma)$, there exists an $\eps > 0$ such that if $\cc$ is a configuration such that $\mu\cc =0$ and $\|SW_2\cc\| < \eps$, then there exists a gauge transformation $g$ such that $g^*\cc \in V_\Sigma$.
\end{enumerate}
\end{Lemma}

\Proof (i) Suppose the statement were not true.  Then we could find a sequence of monopoles $\g_i$ such that $\E_I(\g_i) \to 0$ yet no gauge transformation maps any of the $\g_i$ into $V$.  In particular, since the energies of the configurations $\g_i$ converge, then by \cite[Theorem 5.1.1]{KM}, a subsequence of the $\g_i$ converges in $H^1(I\times\Sigma)$ modulo gauge.  The limiting monopole must have zero energy and therefore belongs to $\fV_I$ modulo gauge.  But this means that for some $i$, a gauge transformation maps $\g_i$ into the neighborhood $V$, a contradiction.

(ii) We have a corresponding energy identity for arbitrary configurations $\cc$ of $\fC(\Sigma)$:
$$\int_{\Sigma} \left(\frac{1}{4}|F_{C}|^2 + |\nabla_C\Upsilon|^2 + \frac{1}{4}\Big(|\Upsilon|^2 + (s/2)\Big)^2 - \frac{s^2}{16}\right) = \|SW_2\cc\|_{L^2(\Sigma)}^2 + \|\mu\cc\|_{L^2(\Sigma)}^2.$$
The proof is now the same as in (i).\End

\begin{Corollary}
  For every gauge invariant neighborhood $V_\Sigma$ of $\fV$ in $\fC^{1/2}(\Sigma)$, there exists an $\eps > 0$ such that if $\g$ is a monopole on $I\times\Sigma$ with $\int_I \|SW_2(\cg(t))\|_{L^2(\Sigma)}^2 < \eps$, then modulo gauge,  we have $\cg(t) \in V_\Sigma$ for all $t \in I$.
\end{Corollary}

\Proof We apply the previous lemma and the embedding $H^1([0,1]\times\Sigma) \hookrightarrow C^0([0,1], H^{1/2}(\Sigma))$.\End

\begin{Lemma}\label{LemmaSmallSF}
  For every $\eps > 0$, there exists an $\eps_0 > 0$ with the following significance.  Let $T \geq 1$ and let $\g$ be a monopole such that $\int_{T-1}^{T+1}\|SW_2(\cg(t))\|_{L^2(\Sigma)} \leq \eps_0$.
  \begin{enumerate}
    \item We have $\|SW_2(\cg(T))\|_{L^2(\Sigma)} \leq \eps.$
    \item If $\cg(T)$ belongs to a coordinate patch $U_\a(\delta)$ for $\delta$ sufficiently small, then $\|\cg(T) - \a\|_{H^s(\Sigma)} \leq C_s\eps$ for all $s \geq 1/2$.
  \end{enumerate}

\end{Lemma}

The above lemma follows easily from the previous lemmas and elliptic bootstrapping, see \cite[Lemma 6.3]{N}.

Given a vortex $\fa$, below are Morse-Bott type inequalities for $CSD^\Sigma_\fa$ in our infinite-dimensional setting.

\begin{Lemma} \label{LemmaMB}
Given a smooth vortex $\a \in \fV$, there exists $\delta > 0$ such the following holds.  If $(C,\Upsilon) \in U_\a(\delta)$ then
    \begin{align}
      |CSD^\Sigma_\a\cc| & \leq \mathrm{const}\cdot \|\cc - \a\|_{H^1(\Sigma)}^2. \label{eq2:MB1}\\
      |CSD^\Sigma_\a\cc|^{1/2} & \leq
      \mathrm{const}\cdot \| SW_2\cc\|_{L^2(\Sigma)}. \label{eq2:MB3}
    \end{align}
\end{Lemma}

\Proof Let $\cco$ and $\cc$ be any two configurations and let $(c,\upsilon) = (C-C_0,\Upsilon-\Upsilon_0)$ be their difference. A simple Taylor expansion of the cubic function $CSD^\Sigma$ shows that it satisfies
\begin{align}
CSD^\Sigma(C_0 + c, \Upsilon_0 + \upsilon) & = CSD^\Sigma(C_0,\Upsilon_0) + ((c,\upsilon), SW_2(C_0,\Upsilon_0)) + \nonumber \\ & \qquad \frac{1}{2}((c,\upsilon),\H_{2,(C_0,\Upsilon_0)} (c,\upsilon)) + \frac{1}{2}(\upsilon,\tilde\rho(c)\upsilon). \label{eq2:CSDdiff}
\end{align}
Letting $\cco$ be a vortex $\a$, then since $SW_2(\a) = 0$, we have from (\ref{eq2:CSDdiff}) that
\begin{align*}
    |CSD_\a^\Sigma(C,\Upsilon)| & \leq \frac{1}{2}\|((c,v), \H_{2,\a}(c,v))\|_{L^2(\Sigma)} + \frac{1}{2}\|(c,v)\|_{L^3(\Sigma)}^3 \\
    & \leq \mathrm{const}\left(\|(c,v)\|_{H^1(\Sigma)}^2 + \|(c,v)\|_{H^{1/2}(\Sigma)\|}\|(c,v)\|_{H^{1}(\Sigma)}^2\right).
\end{align*}
Here, we use that $\H_{2,\a}$ is a first order with smooth coefficients and we use the embedding $H^{1/2}(\Sigma) \hookrightarrow L^4(\Sigma) \subset L^3(\Sigma)$.  The estimate (\ref{eq2:MB1}) now follows from the hypotheses, which implies $\|(c,v)\|_{H^{1/2}(\Sigma)} < \delta$.

The second inequality (\ref{eq2:MB3}) is a standard inequality for Morse-Bott type functions, which one can establish using an infinite-dimensional version of the Morse-Bott lemma (\cite[Chapter 4.5]{Don}). One can see, e.g. by using the same local straightening map analysis of \cite[Lemma 3.4]{N1}, that the Morse-Bott lemma can be performed in a $H^{1/2}(\Sigma)$ neighborhood of the space of vortices.\End

\begin{Remark}\label{RemMB}
  The standard Morse-Bott inequality (in finite dimensions) states that \linebreak $|f(x) - f(a)|^{1/2} \leq c|\nabla_x f|$ holds for all $x$ in a neighborhood of a point $a$ belonging to the critical set of a Morse-Bott function $f$.  In the above, we have been a bit cavalier in our notion of the gradient, since an inner product needs to be specified.  However, since the projection $\Pi_{\K_\a, \J_\cc} : \K_\cc \to \K_\a$ is an isomorphism, uniformly in the $L^2$ norm for $\|\cc - \a\|_{H^{1/2}(\Sigma)}$ sufficiently small, whether we use the usual $L^2$ inner product or the inner product (\ref{eq2:innerprod}) is immaterial.
\end{Remark}

\begin{Definition}
  We say that the chart $U_\a$ is a \textit{Morse-Bott chart} for $\a$ if its closure is contained in a chart of the form $U_\a(\delta)$, with $\delta$ sufficiently small as in Lemma \ref{LemmaMB}.
\end{Definition}

We are interested in configurations which are in standard form with respect to a Morse-Bott coordinate chart.  This is because the Morse-Bott estimates we obtain on these charts allow us to prove the usual exponential decay estimates for Morse-Bott type flows:


\begin{Lemma}
Let $\g$ be a smooth finite energy solution to $SW_3(\g) = 0$ on $Y = [T,\infty) \times \Sigma$ which is in standard form with respect to a Morse-Bott chart $U_\a$ on $Y$.  Then we have the following:
\begin{enumerate}
\item The path $\cg(t)$ converges in $L^2(\Sigma)$ to a vortex $\frak{a}$ as $t \to \infty$ and the temporal component of $\g$ converges in $L^2(\Sigma)$ to zero.
\item The energy of $\g$, or more precisely, the function $CSD_\a^\Sigma(\cg(t))$, decays exponentially as $t \to \infty$.  \end{enumerate}
\end{Lemma}

\Proof Let $CSD^\Sigma_\a(t) = CSD^\Sigma_\a(\cg(t))$.  It is a nonnegative, nonincreasing function of $t$.  We obtain the differential inequality
\begin{align}
  \frac{d}{dt}CSD^\Sigma_\a(t) & = -\left<\Pi_{\K_\a,\J_{\cg(t)}} SW_2(\cg(t)), \Pi_{K_\a,\J_{\cg(t)}}SW_2(\cg(t))\right>_{\a,\cg(t)}\\
  & = -\|SW_2(\cg(t))\|_{L^2(\Sigma)}^2 \\
  & \leq -\mathrm{const}\cdot CSD^\Sigma_\a(t).
\end{align}
In the first line we used (\ref{eq2:PiKJSW2}), in the second line, we used that $\Pi_{\K_{\cg(t)}}\Pi_{\K_\a,J_{\cg(t)}} = \Pi_{\K_{\cg(t)}}$, and in the last line, we used (\ref{eq2:MB3}).
The above inequality implies that
$$CSD^\Sigma_\a(t) \leq c_0e^{-\delta_0t}\cdot CSD^\Sigma_\a(T)$$
for some constants $c_0$ and $\delta_0$ depending on $\a$.  Since the space of vortices is compact modulo gauge however, we can ultimately choose $c_0$ and $\delta_0$ independent of $\a$.

Moreover, we have the following length estimate.  First, we have
\begin{align*}
  \|SW_2(\cg(t))\|_{L^2(\Sigma)} &= \|SW_2(\cg(t))\|^2_{L^2(\Sigma)}\|SW_2(\cg(t))\|^{-1}_{L^2(\Sigma)} \\
  & \leq c\|SW_2(\cg(t))\|^2_{L^2(\Sigma)}CSD_\a^\Sigma(t)^{-1/2}\\
  & = -c\frac{d}{dt}CSD_\a^\Sigma(t)^{1/2},
\end{align*}
Note Remark \ref{RemMB} in passing to the last line.  The above computation makes sense for any non-stationary monopole $\g$, since then $SW_2(\cg(t)) \neq 0$ for every $t$ (otherwise, by unique continuation, we would have $SW_2(\cg(t)) = 0$ for all $t$).  Thus,
\begin{align}
  \|\cg(T_0) - \cg(T_1)\|_{L^2(\Sigma)} & \leq \int_{T_0}^{T_1}\left\|\frac{d}{dt}\cg(t)\right\|_{L^2(\Sigma)}dt \nonumber \\
  &\leq c\int_{T_0}^{T_1}\|SW_2(\cg(t))\|_{L^2(\Sigma)}dt \nonumber \\
  & \leq c'\left(CSD^\Sigma_\a(T_0)^{1/2} - CSD^\Sigma_\a(T_1)^{1/2}\right). \label{eq2:length}
\end{align}
Since $CSD^\Sigma_a(t)$ is decreasing to zero, the $\cg(t)$ form a Cauchy sequence in $L^2(\Sigma)$.  In particular, the path $\cg(t)$ converges to a limit, which must be a vortex.\End

We need two more important facts.  First, we want to show that we can satisfy the hypothesis of the previous lemma, namely that given a monopole with small enough energy, one can always find a coordinate patch about a vortex and a gauge transformation that places the monopole into standard form for all future time with respect to the coordinate patch.  Secondly, we want to show that not only does a monopole in standard form yield a path of configurations convergent in $L^2(\Sigma)$ to a vortex but that the monopole itself on $Y$ converges exponentially fast in all $H^k$ Sobolev norms on $Y$.  This is guaranteed by the following lemma:


\begin{Lemma}\label{LemmaExpDecay}
 There exists an $\eps_0 > 0$ with the following significance.
  \begin{enumerate}
    \item If $\g$ is a monopole such that $\int_T^\infty\|SW_2(\cg(t))\|_{L^2(\Sigma)}^2dt = \eps < \eps_0$, then there exists a Morse-Bott coordinate patch $U_\a$ and a gauge transformation $g$ such that $\g' = g^*\g$ is in standard form with respect to $U_\a$ on $[T,\infty)\times\Sigma$.
    \item There exists a $\delta_0 > 0$ such if $0 < \delta < \delta_0$, then
    \begin{equation}
      \|\g' - \g_{\a'}\|_{H^s([T+t,\infty)\times\Sigma)} \leq C_s C_{\eps}  e^{-\delta t} \label{eq2:energydecay}
    \end{equation}
    for every $s \geq 0$.  Here $\a' = \lim_{t \to \infty}{\check{\g}}'(t)$ is the limiting vortex for $\g'$, $C_s$ is a constant depending on $s$, and $C_{\eps}$ is a constant that can be taken arbitrarily small for $\eps$ sufficiently small.
  \end{enumerate}
\end{Lemma}

\Proof (i) By Lemma \ref{LemmaSmallEnergy}, for $\eps_0$ sufficiently small, we can find a gauge transformation $g$ on $[T,T+1] \times \Sigma$ such that $g^*\g$ is in standard form with respect to some Morse-Bott patch $U_\a$.  The key step is to show that $g$ can be extended to all of $[T,\infty) \times \Sigma$ in such a way that the resulting gauge transformation places $\g$ in standard form for all future time.  However, Lemma \ref{LemmaSmallSF} together with the same arguments as in \cite[Theorem 4.3.1]{MMR} shows that this is the case for $\eps_0$ sufficiently small.

(ii) By (i) and Lemma \ref{LemmaSmallSF}, we know that $\sup_{t \geq T+1}\|\cg'(t) - \a'\|_{H^k(\Sigma)} \leq C_kC_\eps$.  Now standard exponential decay arguments, e.g. \cite[Lemma 5.4.1]{MMR}\footnote{Note this lemma is a more general statement than we need, since in our Morse-Bott situation, the center manifold is simply the critical manifold.} yields the desired conclusion for $s=0$. For $s > 0$, we use the fact that one can bootstrap elliptic estimates in the standard form gauge so that $L^2$ exponential decay gives us $H^s$ decay on the cylinder.  The arguments are formally similar to those of \cite[Lemma 3.3.2]{MMR}.  We omit the details.  Note that we can take $\delta_0$ independent of $\a$ since the vortex moduli space $\V$ is compact.\End


\section{The Finite Energy Moduli Space}


In this section, we use the results developed in the previous section to prove our main results concerning the space of finite energy monopoles on a $3$-manifold with cylindrical ends.  We first establish the basic setup (i.e. the appropriate configuration spaces and gauge groups) based on our previous analysis and then divide our task into first studying the case when $Y$ is a semi-infinite cylinder $[0,\infty) \times \Sigma$ and then proceeding to the general case.

From Lemma \ref{LemmaExpDecay}, we see that modulo gauge, any finite energy monopole converges exponentially to a vortex.  This result depends crucially on the Morse-Bott framework\footnote{For comparison, in \cite{MMR}, one does not always get exponential decay in the instanton case due to Morse-Bott degenerate critical points.} of the previous section and it yields for us the following two bits of information. First, it tells us that the right choice of function spaces to consider on the cylinder are the exponentially weighted Sobolev spaces.  Second, it suggests that the topology of our monopole spaces is related to the topology of the vortex moduli spaces at infinity.  Our main theorem  of this section is Theorem \ref{ThmFinEnergy}.  We also prove Theorems \ref{ThmSmallEnergy} and \ref{ThmLagVort} to describe other moduli spaces one might wish to obtain from boundary values of monopoles.

When $Y$ is a manifold with cylindrical ends (and without boundary), every monopole on it is gauge equivalent to a smooth monopole.  However, our preliminary analysis requires us to work in the general Hilbert space setting (so that as usual, Banach space methods can be employed).  On the noncompact space $Y$, we must a priori work with local Sobolev spaces $H^s_\loc(Y)$, that is, the topological vector space of functions on $Y$ that belong to $H^s(K)$ for every compact domain $K \subset Y$.  We let $\fC^s_\loc(Y)$ denote the $H^s_\loc(Y)$ completion of the smooth configuration space on $Y$.  Then the space of all finite energy monopoles in $\fC^s_\loc(Y)$ is given by \label{p:fMfinenergy}
$$\fM^s = \fM^s_\eta(Y) = \{\g \in \fC^s_\loc(Y) : SW_3(\g) = \eta, \E(\g) < \infty\}.$$
As before, we assume $\eta$ is of the special form $\eta = -\eta^0 dt$ along each end, with $\eta^0$ time-translation invariant. We topologize this space in the $H^s_\loc(Y)$ topology and also by requiring that the energy be a continuous function. Likewise, for any $E > 0$, we can define the space
$$\fM^s_E = \{\g \in \fM^s : \E(\g) < E\}.$$
of $H^s_\loc(Y)$ monopoles that have energy less than $E$.  Here, $\E$ is the analytic (equivalently, topological) energy of \cite{KM}, which when restricted to the cylindrical ends of $Y$, is given by the expression (\ref{eq:energy}).

These spaces, being merely the spaces which a priori contains all the monopoles of interest, are much too large to be of use.  Of course, as we have mentioned, we can always find a gauge in which a finite energy configuration decays exponentially in every Sobolev norm at infinity.  So for any $\delta \in \R$ and nonnegative integer $s \geq 0$, define $H^{s;\delta}(Y)$ to be the closure of $C^\infty_0(Y)$ in the norm \label{p:Hsdelta}
$$\|f\|_{H^{s;\delta}(Y)} = \|e^{\delta t}f\|_{H^s(Y)}.$$
Thus, for $\delta > 0$, the weight $e^{\delta t}$ forces exponential decay of our functions; for $\delta < 0$, we allow exponential growth.  Using this topology, we can topologize the space $\T = \Omega^1(Y; i\R) \oplus \Gamma(\S)$ (the tangent space to the smooth configuration space on $Y$ when $Y$ was compact) in the $H^{s; \delta}(Y)$ topology to obtain $\T^{s;\delta}$.  For $\delta > 0$, we can then define the corresponding space
$$\fC^{s;\delta}(Y) = \{\g : \g - \g_\a \in \T^{s;\delta} \textrm{ for some }\a \in \fV^s\}$$
of configurations that decay exponentially to some $H^s(\Sigma)$ vortex $\a \in \fV^s := H^s\fV(\Sigma)$. (As we have been doing consistently, to avoid notational clutter, we suppress the dependence of $\fV(\Sigma)$ on the connected components of $\Sigma$, the $\spinc$ structure, and the perturbation. The dependence is given explicitly by Lemma \ref{LemmaVortices}.)
In particular, if $s \geq 2$, all configurations in $\fC^{s;\delta}(Y)$ are pointwise bounded.  From now on, we will assume $s$ is an integer and $s \geq 2$ unless otherwise stated.  We give $\fC^{s;\delta}(Y)$ the topology of $\T^{s;\delta} \times \fV^s$ in the obvious way.  In particular, observe that $\fC^{s;\delta}(Y)$ is a Hilbert manifold.  
Given $\g \in \fM^s$, the gauge transformation which sends $\g$ to an element of $\fC^{s;\delta}(Y)$, being only required to satisfy a condition at infinity, can be taken to be identically one near $\Sigma = \partial Y$. It follows that to study the space $\fM^s$ and its boundary values on $\fC^{s-1/2}(\Sigma)$, it suffices to study the space
$$\fM^{s;\delta} = \fM^s \cap \fC^{s;\delta}(Y),$$
for $\delta > 0$ small.

We now consider the above setup modulo all gauge transformations.  When $\delta > 0$, the exponential decay of configurations allows multiplication to be possible and we can define an exponentially weighted gauge group accordingly.  Namely, we define $\G^{s+1;\delta}(Y)$ to be the Hilbert Lie group of gauge transformations such that $g$ differs from a constant gauge transformation by an element of $H^{s+1; \delta}(Y)$. This group acts smoothly on $\fC^{s;\delta}(Y)$ and we can form the quotient space \label{p:fCdelta}
$$\fB^{s;\delta}(Y) = \fC^{s;\delta}(Y)/\G^{s+1;\delta}(Y).$$
This quotient is a smooth Hilbert manifold away from the reducible configurations, which we can ignore when studying the monopole moduli space since none of our vortices are reducible by the Morse-Bott hypothesis.  Moreover, it is diffeomorphic to a Hilbert bundle over the space
$$\V_Y(\Sigma) := \fV^s(\Sigma)/\left(\G^{s+1;\delta}(Y)|_\Sigma\right),$$
which is a covering space of $\V(\Sigma)$ whose fiber is the image of $H^1(Y; \Z)$ inside $H^1(\Sigma; \Z)$ under the natural restriction map (i.e. the component group of those gauge transformations on $\Sigma$ that extend to $Y$).

Let \label{p:romanM}
\begin{equation}
  M^s = M^s_\eta(Y) = \fM^{s;\delta}_\eta(Y)/\G^{s+1;\delta}(Y) \subset \fB^{s;\delta}(Y) \label{eq2:MsY}
\end{equation}
denote the moduli space of gauge equivalence classes of monopoles on $Y$ which decay exponentially to a vortex.  Note that by our exponential decay results, $M^s$, topologized as a subspace of $\fB^{s;\delta}(Y)$, is also (topologically) the quotient space of $\fM^s$ by the group of $H^{s+1}_\loc(Y)$ gauge transformations on $Y$.  (Here, it is key that $\fM^s$ is topologized with the energy functional.) Observe that the definition of $M^s$ is independent of $\delta$ for $\delta > 0$ sufficiently small as a consequence of Lemma \ref{LemmaExpDecay}.

\subsection{The Semi-Infinite Cylinder}

We first specialize to $Y = [0,\infty) \times \Sigma$ to isolate the analysis on the cylindrical end.  First, we obtain tangent space decompositions of our configuration space arising from the infinitesimal gauge action, as in \cite{KM} and \cite{N1}, but on weighted spaces.  These decompositions are important for studying the linearization of the gauge-fixed Seiberg-Witten equations.  Thus, for $\c \in \fC^{s;\delta'}(Y)$ with $\delta' > 0$, we can define the operators
\begin{align*}
  \bd_\c: H^{s+1;\delta}(Y;i\R) & \to \T^{s;\delta}\\
  \xi & \mapsto (-d\xi,\xi\Psi)\\
  \bd_\c^*: \T^{s;\delta} & \to \T^{s-1;\delta}\\
  (b,\psi) & \mapsto -d^*b + i\Re(i\Psi,\psi).
\end{align*}
for $\delta \leq \delta'$.  We then obtain the subspaces
\begin{align*}
  \J_\c^{s;\delta} &= \im \bd_\c\\
  \K_\c^{s;\delta} &= \ker \bd_\c^*
\end{align*}
of $\T_\c^{s;\delta}$.  With the presence of a boundary $\partial Y = \Sigma$, we can also supplement the above subspaces with boundary conditions:
\begin{align}
  \J^{s;\delta}_{\c,t} &= \{(-d\xi,\xi\Psi) : \xi \in H^{s+1;\delta}(Y;i\R),\; \xi|_\Sigma = 0\},\\
  \K^{s;\delta}_{\c,n} &= \{(b,\psi) \in \K^{s;\delta}_\c : *b|_\Sigma = 0\},
\end{align}
By standard Fredholm theory on weighted spaces (see \cite{LM}), we have a weighted decomposition
\begin{equation}
  \T^{s;\delta} = \J^{s;\delta}_\c \oplus \K^{s;\delta}_{\c,n} \label{eq2:0-T=JKn}
\end{equation}
for $\c$ irreducible.  This is summarized in Lemma \ref{LemmaTw}.

Define the operator
\begin{align}
  \partial_\infty: \fB^{s;\delta}(Y) & \to \V(\Sigma) \nonumber \\
  [\g] & \mapsto \lim_{t \to \infty}[\cg(t)] \label{eq2:partialB}
\end{align}
mapping a configuration to its limiting vortex at infinity. Given any irreducible $\g \in \fC^s(Y)$, from (\ref{eq2:0-T=JKn}), we have that the tangent space to $[\g]$ of $\fB^{s;\delta}(Y)$ can be identified with
\begin{equation}
  T_{[\g]}\fB^{s;\delta}(Y) \cong \K^{s;\delta}_{\g,n} \cap T_{\partial_\infty[\g]}\V. \label{eq2:TfB}
\end{equation}
The map (\ref{eq2:partialB}) restricts to a map
$$\partial_\infty: M^s(Y) \to \V$$
mapping a monopole to its asymptotic vortex on $\Sigma$.  Given a vortex $[\a] \in \V$, we can define
\begin{equation}
  M^s(\a) = \{[\g] \in M^s(Y) : \partial_\infty[\g] = [\a]\}
\end{equation}
the moduli space of monopoles that converge to $[\a]$.

We are now in the position to state our main result.  We have the (tangential) restriction map
\begin{align*}
  r_\Sigma: \fC^{s;\delta}(Y) & \to \fC^{s-1/2}(\Sigma) \\
  (B,\Psi) & \mapsto (B,\Psi)|_\Sigma,
\end{align*}
where $\Sigma = \{0\}\times\Sigma$.
Letting
\begin{equation}
  \fB^s(\Sigma) = \fC^s(\Sigma)/\G^{s+1}(\Sigma)
\end{equation}
denote the quotient configuration space on $\Sigma$, the restriction map $r_\Sigma$ descends to the quotient space:
\begin{equation}
  r_\Sigma: \fB^{s;\delta}(Y) \to \fB^{s-1/2}(\Sigma).
\end{equation}
Let
\begin{equation}
  \fB_{\eta^0}^s(\Sigma) = \mu^{-1}(\eta^0)/\G^{s+1}(\Sigma) \subset \fB^s(\Sigma)
\end{equation}
denote the symplectically reduced space associated to the $\eta^0$ level set of moment map $\mu$, where $\eta^0 \in \Omega^0(\Sigma; i\R)$.  We have the following theorem, which geometrically, is the statement that $M^s([0,\infty) \times \Sigma)$ is the (infinite-dimensional) stable manifold to the space of vortices on $\Sigma$ under the Seiberg-Witten flow.

\begin{Theorem}\label{ThmFinEnergy} (Finite Energy Moduli Space) Let $Y = [0,\infty) \times \Sigma$, fix a $\spinc$ structure $\s$ on $\Sigma$, and let $s \geq 2$ be an integer.  Let $d = \frac{1}{2}\left<c_1(\s),\Sigma\right>$ and pick $\eta^0$ such that $\frac{i}{\pi}\int_\Sigma \eta^0 \neq d$. Then the following holds:
  \begin{enumerate}
    \item The moduli space $M^s(Y) = M^s_{\eta^0}(Y)$ is naturally a smooth Hilbert submanifold\footnote{From now on, we will always regard $M^s(Y)$ as endowed with this topology. As mentioned, it is homeomorphic to the quotient of $\fM^s$ by $\G^{s+1}_\loc(Y)$, but this latter space does not come with an a priori smooth manifold structure.} of $\fB^{s;\delta}(Y)$, for $\delta > 0$ sufficiently small.
    \item The map $r_\Sigma: M^s(Y) \to \fB^{s-1/2}(\Sigma)$ is a diffeomorphism onto its image, which is a coisotropic submanifold of the symplectically reduced space $\fB^{s-1/2}_{\eta^0}(\Sigma)$.  Given any $[\g] \in M^s(Y)$, the annihilator of the coisotropic space $r_\Sigma(T_{[\g]}M^s(Y))$ is the space $r_\Sigma\big(T_{[\g]}M^s(\partial_\infty[\g])\big)$.
    \item Both $M^s(Y)$ and $r_\Sigma(M^s(Y))$ are complete.
  \end{enumerate}
\end{Theorem}

In regarding $M^s(Y)$ as the stable manifold to the space of vortices at infinity, we see that it is the union of the $M^s(\a)$, each of which is the stable manifold to $[\a] \in \V$, as $[\a]$ varies over the symplectic set of critical points $\V$.  This geometric picture clarifies the symplectic nature of (ii) in the above.

Because of the infinite-dimensional nature of the objects involved, the proof of the above theorem requires some care.  We first prove a few lemmas.  The first lemma below extends tangent space decompositions to weighted spaces, which is needed in understanding transversality for the Seiberg-Witten map $SW_3$ as a section of the relevant Hilbert bundle. Given a configuration $\g \in \fC^{s;\delta}(Y)$, the linearization of the Seiberg-Witten map $SW_3: \fC^{s;\delta}(Y) \to \T^{s-1;\delta}$ produces for us a first order formally self-adjoint operator, the \textit{Hessian}:
\begin{equation}
  \H_\g: \T^{s;\delta} \to \T^{s-1;\delta}.
\end{equation}

\begin{Lemma}\label{LemmaTw}
  Let $s \geq 2$ and $\g \in \fC^{s;\delta'}(Y)$ where $\delta' > 0$.  Then for $\delta > 0$ sufficiently small, the following hold:
  \begin{enumerate}
  \item We have the following decomposition for $1 \leq s' \leq s$:
    \begin{align}
        \T^{s';\pm\delta} &= \J^{s';\pm\delta}_{\g,t} \oplus \K^{s';\pm\delta}_\g. \label{eq2:T=JKw}
    \end{align}
    If $\g$ is not reducible, then we also have
    \begin{align}
        \T^{s';\pm\delta} &= \J^{s';\pm\delta}_{\g} \oplus \K^{s';\pm\delta}_{\g,n}. \label{eq2:T=JKwn}
    \end{align}
  \item Let $SW_3(\g) = 0$. Then the Hessian operator $\H_\g: \T^{s;\pm\delta} \to \K^{s-1;\pm\delta}$ is surjective.
  \end{enumerate}
\end{Lemma}

\Proof Using the Fredholm theory for elliptic operators on weighted spaces of \cite{LM}, the proof of this lemma proceeds mutatis mutandis as in Lemmas 3.4, 3.16, and 4.1 of \cite{N1}, since the elliptic methods there adapt to weighted spaces for weights on the complement of a discrete set.\End

Because of (ii) above, we have that the space of monopoles on $Y \times \Sigma$ is transversally cut out by $SW_3$.  We now investigate the symplectic aspects of the boundary values of the space of monopoles.

To take into account gauge-fixing and obtain an elliptic operator from the Hessian, we proceed analogously to \cite{KM} and define the \textit{extended Hessian}
  \begin{equation}
   \wH_\g := \begin{pmatrix}
      \H_\g & \bd_\g \\
      \bd_\g^* & 0
    \end{pmatrix}: \T^{s;\delta} \oplus H^{s;\delta}(Y;i\R) \to \T^{s-1;\delta} \oplus H^{s;\delta}(Y;i\R).
  \end{equation}
That is, $\wH_\g$ is obtained from $\H_\g$ by taking into the account the gauge action of the exponentially decaying gauge group $\G^{s+1;\delta}(Y)$.

Thus, define the following augmented space
$$\tT^{s;\delta} = \T^{s;\delta} \oplus H^{s;\delta}(Y;i\R)$$
for $\delta \in \R$.  In general, for $\g \in \fC^{\delta'}(Y)$ and $\delta \leq \delta'$, we have the first order formally self-adjoint elliptic operator
$$\hH_\g: \tT^{s;\delta} \to \tT^{s-1;\delta}.$$
In \cite{N1}, a method known as the ``invertible double" (see \cite{BBLZ}) is used as a fundamental tool in showing that the space of boundary values of kernel of the operator\footnote{Or rather, the operator referred to as the augmented Hessian there.} $\wH_\g$ on a compact $3$-manifold yields a Lagrangian subspace of the boundary data space.  Here, the same methods can be used, only now we have a slightly different situation due to the weights. Nevertheless, this invertible double technique is what allows us to obtain symplectic information for the boundary data of the kernel of the extended Hessian in the cylindrical case.

As in \cite{N1}, we first observe that
$$\tT_\Sigma := \T_\Sigma \oplus \Omega^0(\Sigma; i\R) \oplus \Omega^0(\Sigma; i\R)$$
is the full boundary value space of $\tT$.  Here, we have a full restriction map $r$ given by
\begin{align}
  r: \tT & \to \tT_\Sigma \nonumber \\
  (b,\psi,\alpha) & \mapsto (b|_\Sigma, \psi|_\Sigma, -\partial_t \llcorner b, \alpha|_\Sigma), \label{rmap}
\end{align}
i.e., the first two components of $r$ is just the tangential restriction map $r_\Sigma$, the next component is the outward normal component of $b$, and the final component is the boundary value of $0$-form $\alpha$.   Extending to Sobolev spaces, we have $r: \tT^{s;\delta} \to \tT_\Sigma^{s-1/2}$ for $s > 1/2$.  We also have the complex structure $\tJ: \tT_\Sigma \to \tT_\Sigma$, given by
$$\tJ = J \oplus \begin{pmatrix}0 & 1 \\ -1 & 0 \end{pmatrix}$$
which extends the complex structure $J$ on $\T_\Sigma$ and which is compatible with the product symplectic form
$$\tomega((b,\phi,\alpha_1,\alpha_2), (a,\psi,\beta_1,\beta_2)) = \omega((b,\phi),(a,\psi)) - \int_\Sigma (\alpha_1\beta_0 - \alpha_0\beta_1)$$
on $\tT_\Sigma$.  As in \cite{N1}, we obtain symplectic data on $\T_\Sigma$ from symplectic data on $\tT_\Sigma$ via symplectic reduction with respect to the coisotropic space $\T_\Sigma \oplus \Omega^0(\Sigma; i\R) \oplus 0$.  Hence, we first study symplectic data on $\tT_\Sigma$, where we can use elliptic methods, in particular, the invertible double method.

\begin{Lemma}\label{LemmaWID}(Weighted Invertible Double)
  Let $s \geq 2$ and let $\delta \neq 0$ be sufficiently small.  Let $\g \in \M^{s;\delta}$.  Define
  $$\tT^{s;\pm \delta} \oplus_{\tJ} \tT^{s;\pm\delta} = \{(x,y) \in \tT^{s;\pm \delta} \oplus \tT^{s;\pm\delta} : r(x) = \tilde J r(y)\},$$
Then we have the following:
  \begin{enumerate}
    \item The ``doubled operator''
    $$\wH_\g \oplus \wH_\g : \tT^{s; \delta} \oplus_{\tJ} \tT^{s;-\delta} \to \tT^{s-1; \delta} \oplus \tT^{s-1;-\delta}$$
    is an isomorphism.
    \item The space $r(\ker \wH_\g|_{\tT^{s;\delta}})$ is an isotropic subspace of $\tT^{s-1/2}_\Sigma$.  Its symplectic annihilator is the coisotropic subspace $r(\ker \wH_\g|_{\tT^{s;-\delta}})$.
  \end{enumerate}
\end{Lemma}

\Proof (i) One can easily construct a parametrix for the double using the methods of \cite{APS}.  This shows that the double is Fredholm.  Here $\delta \neq 0$ small is needed because of our Morse-Bott situation at infinity.  To see that the double is injective, if $u = (u_+,u_-) \in \tT^{s; \delta} \oplus_{\tJ} \tT^{s;-\delta}$ belongs to the kernel of the double, then
$$0 = (u_+, \wH_\g u_-)_{L^2(Y)} - (\wH_\g u_+, u_-)_{L^2(Y)}  = -\int_\Sigma \left(r(u_+), \tJ r(u_-)\right).$$
The second equality is Green's formula (i.e., an integration by parts formula) for $\wH_\g$, where $\Sigma = \partial Y$.  This formula is justified since $u_+$ decays exponentially while $u^-$ is at most bounded since $\delta$ is small (see (\ref{eq2:kernelspace})), so that there is no contribution from infinity.  On the other hand, since $r(u_+) = \tJ r(u_-)$, we conclude that
$$\int_\Sigma |u_+|^2 = \int_\Sigma |u_-|^2 = 0.$$
Thus $u = 0$ and so the double is injective.  Integration by parts and the same argument shows that the orthogonal complement of the range of the double is zero.  Thus, the double is invertible.

(ii) Green's formula above shows that $r(\ker \wH_\g|_{\tT^{s;\delta}})$ is isotropic and that it annihilates $r(\ker \wH_\g|_{\tT^{s;-\delta}})$.  It remains to show that the annihilator of $r(\ker \wH_\g|_{\tT^{s;\delta}})$ is precisely $r(\ker \wH_\g|_{\tT^{s;-\delta}})$, for which it suffices to show that $r(\ker \wH_\g|_{\tT^{s;\delta}})$ and $\tJ r(\ker \wH_\g|_{\tT^{s;-\delta}})$ are (orthogonal) complements.  This however follows from (i) and the same method of proof of \cite[Proposition 5.12]{BBLZ}.\End


\begin{Lemma}\label{LemmaZero}
  Let $\g_\a$ be the configuration on $Y$ corresponding to the constant path in $\fC(\Sigma)$ identically equal to a configuration $\a$ (not necessarily a vortex).
\begin{enumerate}
  \item We can write
  \begin{equation}
    \wH_{\g_\a} = \tJ\left(\frac{d}{dt} + \barB_{\a}\right)
  \end{equation}
  where $\barB_{\a}: \tT_\Sigma \to \tT_\Sigma$ is a time-independent first order self-adjoint operator given by
  \begin{equation}
    \barB_\a = \begin{pmatrix}
  \H_{2,\a} & \bd_\a & -J \bd_\a \\
  \bd_\a^* & 0 & 0 \\
  \bd_\a^* J & 0 & 0
\end{pmatrix}
  \end{equation}
  \item We have $\tJ \barB_{\a} = -\barB_{\a} \tJ $.
  \item If $\a = (C,\Upsilon)$ where $\Upsilon = (\Upsilon^+,0)$, then $\barB_{\a}$ is complex-linear with respect to the complex structure $I: \tT_\Sigma \to \tT_\Sigma$ given by
      $$I = (\check{*}, i) \oplus \begin{pmatrix}0 & 1 \\ -1 & 0\end{pmatrix}.$$
      If instead $\Upsilon = (0,\Upsilon^-)$, the analogous statement is true with $i$ replaced with $-i$ in the above.
  \item If $\a$ is a vortex, then
  \begin{equation}
      \ker \barB_\a = \{(c,\upsilon,0,0) \in \tT_\Sigma: (c,v) \in T_\a\fV, \;
      \bd_\a^*(c,\upsilon)=0\}  \label{eq2:kerB}
  \end{equation}
    is isomorphic to the tangent space to the vortex moduli space $\V$ at $[\a]$.

\end{enumerate}
\end{Lemma}

\Proof (i) This is a straightforward computation. Observe that the operators appearing in the first column of $\bB_\a$ correspond to the linearization of $SW_2$, gauge fixing, and the moment map, respectively.  (ii, iii) Using the fact that $J$ and $(\check{*}, i)$ anti-commute and commute with $\H_{2,\a}$, respectively, (with the appropriate assumption on $\a$ in the latter case) this is a straightforward computation.

(iv) Suppose $(c,\beta,\alpha,\upsilon) \in \ker \barB_{\a}$. Then
$$\H_{2,\a}(c,\upsilon) - J\bd_\a\beta + \bd_\a\alpha = 0.$$
All three terms in the above however are orthogonal to each other, since the tangent space to the gauge group is isotropic and since Proposition \ref{PropMomentMap}(iv) holds.  It follows that $J\bd_\a\beta = \bd_\a\alpha = 0$, whence $\beta = \alpha = 0$ since $\a$ is not reducible.  We now have the equality (\ref{eq2:kerB}), since being annihilated by the first column of $\barB_\a$ expresses being a gauge-fixed element of the kernel of the linearized vortex equations.\End

Note that the special algebraic structure of $\bB_\a$ comes from its symplectic origins.  Indeed, $\bB_\a$ possesses an anti-commuting complex structure by virtue of arising from the tangential part of the Dirac operator $\bH_{\g_\a}$.  Furthermore, $\bB_\a$ possesses a commuting complex structure by virtue of its association to the vortex equations, which can be interpreted as the zero level set of a moment map on a K\"ahler configuration space \cite{GP} (when $\a = (C,\Upsilon^+)$ is not a vortex, the vanishing of $\Upsilon^-$ is enough to preserve the complex-linearity of $I$). \\


\noindent \textbf{Proof of Theorem \ref{ThmFinEnergy}:} (i) We prove that $\M^{s;\delta}$ is a Hilbert manifold by showing that it is the zero set of a section of a Hilbert bundle that is transverse to zero.  We have the exponentially decaying space $\K^{s-1;\delta}_\g$ for every configuration $\g$ on $Y$.  Just as in \cite[Proposition 3.5]{N1}, since $\K^{s-1;\delta}_\g$ varies continuously with $\g$, we may form the bundle $\K^{s-1;\delta}(Y) \to \fC^{s;\delta}(Y)$ whose fiber over every $\g \in \fC^{s;\delta}(Y)$ is the Hilbert space $\K^{s-1;\delta}(Y)$. We can interpret $SW_3$ as a section of the bundle $\K^{s-1;\delta}(Y)$, i.e.,
\begin{align}
  SW_3: \fC^{s;\delta}(Y) \to \K^{s-1;\delta}(Y). \label{eq2:SW3trans}
\end{align}
Note that the range of $SW_3$ really is contained in the exponentially decaying space $\K^{s;\delta}(Y)$.  Indeed, for any constant vortex $\g_\a$ induced from $\a \in \fV^s$ and any $x \in \T^{s;\delta}$, we have
$$SW_3(\g_\a + x) = \H_{\g_\a}x + x \sharp x,$$
where $\sharp$ denotes a bilinear pointwise multiplication map.  Since $\delta > 0$ and $s \geq 2$, multiplication is bounded on $\T^{s;\delta}$ and so in particular, $x \sharp x \in \T^{s-1;\delta}$.  Lemma \ref{LemmaTw}(ii) implies (\ref{eq2:SW3trans}) is transverse to the zero section, whence $\fM^{s;\delta} = SW_3^{-1}(0)$ is a smooth Hilbert submanifold of $\fC^{s;\delta}(Y)$. Since there are no reducibles, $\G^{s+1;\delta}(Y)$ acts freely, and so $M^s \cong \fM^{s;\delta}/\G^{s+1;\delta}(Y)$ has the structure of a smooth Hilbert submanifold of $\fB^{s;\delta}(Y)$.

(iii) We have $\fM^{s;\delta} \cong M^s \times \G^{s+1;\delta}(Y)$ (cf. \cite[Chapter 9.3]{KM}), so that since $\fM^{s;\delta}$ is complete, so is $M^s$. To show that $r_\Sigma(M^s)$ is complete, we have to show that any sequence in $r_\Sigma(M^s)$ which forms a Cauchy sequence in $\fB^{s-1/2}(\Sigma)$ converges to an element of $r_\Sigma(M^s)$.  Since $s - 1/2 \geq 1/2$, if a sequence converges in $H^{s-1/2}(\Sigma)$, its values under $CSD^\Sigma$ converge, since $CSD^\Sigma$ is $H^{1/2}(\Sigma)$ continuous.  Thus, it follows that the limiting configuration has finite energy, and the limiting trajectory it determines on the cylinder is the limit of the sequence of trajectories.  Thus, the limit corresponds to a finite energy monopole, and hence $r_\Sigma(M^s)$ is complete.

(ii) For the first part of (ii), similar unique continuation arguments as made in the proof of the main theorem of \cite{N1} imply the injectivity of $r_\Sigma$ and that it is an immersion.  To show then that $r_\Sigma$ is a global embedding, we use similar arguments as made in the proof of \cite[Theorem 4.13]{N1}.  Namely, it suffices to show that if $r_\Sigma(\g_i)$ forms a Cauchy sequence in $\fB^{s-1/2}(\Sigma)$, then the $\g_i$ form a Cauchy sequence in $M^s$.  However, this follows from our preceding analysis.  Namely, we have that the energy of the $\g_i$ converge.  On compact cylinders, the $\g_i$ converge in $H^s(I\times\Sigma)$ by \cite[Lemma 4.11]{N1}, and at infinity, we have convergence in a neighborhood of infinity since energy controls exponential decay, i.e., we have equation (\ref{eq2:energydecay}).  Because of the way $M^s$ is topologized, this gives us convergence of $\g_i$ in $M^s$.

It remains to prove the more interesting second part of (ii).  Let $\H_\g^{s;\pm \delta}$ and $\barH_\g^{s;\pm\delta}$ be the Hessian and extended Hessian operators with domains $\T^{s;\pm\delta}$ and $\tT^{s;\pm\delta}$, respectively. Observe that given $[\g] \in M^s$, then $r_\Sigma (T_{[\g]}M)$ can be regarded as the symplectic reduction\footnote{See \cite{N1} for further reading.} of $r_\Sigma(T_\g\fM^{s;\delta})$ with respect to the coisotropic subspace $T_\g\mu^{-1}(\eta^0)$ of $T^{s-1/2}_\Sigma$.  We have the following claim:\\

\noindent \textit{Claim: } The space $r_\Sigma(T_\g\fM^{s;\delta})$ is a coisotropic subspace of $\T^{s-1/2}_\Sigma$ with annihilator $r_\Sigma(T_\g\fM^{s;\delta} \cap \T^{s;\delta}).$\\

We will prove this claim, which is equivalent to second assertion of (ii) via the previous observation concerning symplectic reduction. To prove the first part of the claim, we proceed as follows.  First of all, we have
\begin{align}
  r_\Sigma(T_\g\fM^{s;\delta}) 
  &= r_\Sigma(\ker \barH_\g |_{T_\g\fC^{s;\delta}(Y)}) \label{eq2:kerHonTC}
\end{align}
which follows from Lemma \ref{LemmaTw} and the fact that $\J^{s;\delta}_{\g,t}$ has zero restriction to the boundary.  Let $\pi_{SR}: \tT_\Sigma^{s-1/2} \to \T_\Sigma^{s-1/2}$ denote the symplectic reduction induced by the coisotropic subspace $W := \T^{s-1/2}_\Sigma \oplus \Omega^0(\Sigma; i\R) \oplus 0$, that is, $\pi_{SR}(x)$ is coordinate projection onto $\T_\Sigma^{s-1/2}$ if $x \in W$ and $\pi_{SR}(x) = 0$ otherwise.  We will show that
\begin{equation}
  r_\Sigma(\ker \barH_\g |_{T_\g\fC^{s;\delta}(Y)}) = \pi_{SR}r(\ker \barH_\g^{s;-\delta/2}), \label{eq2:kerH-SR}
\end{equation}
which together with Lemma \ref{LemmaWID} and (\ref{eq2:kerHonTC}) will show that $r_\Sigma(T_\g\fM^{s;\delta})$ is coisotropic.

Let $\a = \lim_{t \to \infty}\cg(t)$.  Then we can write
\begin{equation}
  \barH_{\g} = \barH_{\g_\a} + R
\end{equation}
where $\barH_{\g_\a}$ is time-independent and where $R$ is a zeroth order operator whose coefficients belong to $H^{s;\delta}(Y)$.  From this, we have
\begin{equation}
  \ker \barH_{\g}^{s;-\delta/2} = \{x \in \tT^{s;\delta/2} + \ker \barH_{\g_\a}^{s;-\delta/2}: \barH_{\g}x = 0\}. \label{eq2:kerH's}
\end{equation}
Indeed, if $x \in \tT^{s;-\delta/2}$ and $\barH_{\g}x = 0$, then $\barH_{\g_\a}x = -Rx \in \T^{s;\delta/2}$.  The operator $\barH_{\g_\a}: \tT^{s+1,\delta/2} \to \tT^{s,\delta/2}$ is surjective (since no boundary conditions are specified), and hence we see that $x$ differs from an element of $\tT^{s+1,\delta/2}$ by an element of $\ker \barH_{\g_\a}^{s;-\delta/2}$.  

Let $Z_{\barB_\a} \subset \ker \barH_{\g_\a}^{s;-\delta/2}$ denote the time-translation invariant elements given by the zero eigenspace of $\barB_\a$. Then for $\delta$ sufficiently small,
\begin{equation}
  \ker \barH_{\g}^{s;-\delta/2} \subset \tT^{s;\delta/2} + Z_{\barB_\a}. \label{eq2:kernelspace}
\end{equation}
Since $Z_{\barB_\a} \subset T_{\g_\a}\fV_{[0,\infty)}$ by Lemma \ref{LemmaZero}, equations (\ref{eq2:kerH's}) and (\ref{eq2:kernelspace}) imply
\begin{equation}
  \ker \barH_{\g}^{s;-\delta/2} = \ker \barH_{\g}|_{T_\g\fC^{s;\delta/2} \oplus H^{s;\delta/2}(Y)}. \label{eq2:kerHT}
\end{equation}
Summarizing the above, we have shown that the only elements of $\barH_{\g}^{s;-\delta/2}$ that do not exponentially decay are those that have a nonzero contribution from  $Z_{\barB_\a} \subset T_\a\fV$.  Because of the direct sum decomposition (\ref{eq2:T=JKw}), elements of $\ker \barH_{\g}^{s;-\delta/2}$ whose restriction under $r$ lie inside the coisotropic space $\T^{s-1/2}_\Sigma \oplus \Omega^0(\Sigma; i\R) \oplus 0$ have vanishing $H^{s;-\delta/2}(Y)$ component, and thus belong to $T_\g\fC^{s;\delta/2}$. (This is exactly the same type of analysis carried out in the symplectic aspects of \cite[Section 3.3]{N1}).  This observation together with (\ref{eq2:kerHT}) implies
$$\pi_{SR}r(\ker \barH_\g^{s;-\delta/2}) = \pi_{SR}r\ker (\barH_{\g}|_{T_\g\fC^{s;\delta/2}}).$$
But we have
$$\pi_{SR}r\ker (\barH_{\g}|_{T_\g\fC^{s;\delta/2}}) = r_\Sigma\ker (\barH_{\g}|_{T_\g\fC^{s;\delta/2}}),$$
and so (\ref{eq2:kerH-SR}) follows from the above two equations.  This finishes the first part of the claim.

The second part of the claim is now a simple consequence of Lemma \ref{LemmaWID}(ii) and the preceding analysis. Namely, we have that the annihilator of $r_\Sigma(T_\g\fM^{s;\delta}) = \pi_{SR}r(\ker \barH_\g^{s;-\delta/2})$ is given by
\begin{align*}
  \pi_{SR}r(\ker \barH_\g^{s;\delta/2}) &=  \pi_{SR}r(\ker \barH_\g|_{\T^{s;\delta/2}})\\
  &= r_\Sigma(T_\g\fM^{s;\delta} \cap T^{s;\delta/2}).\label{eq2:stablespace}
\end{align*}
The claim now follows from the fact that $T_\g\fM^{s;\delta} \cap T^{s;\delta/2}$ modulo gauge is precisely $T_{[\g]}M^s({\partial_\infty [\g]})$.\End

The next two results we state concern two natural additional moduli spaces one may consider: those monopoles with small energy and those monopoles whose limiting value belong to some Lagrangian submanifold of $\V(\Sigma)$.  For the first of these, we consider the moduli space
$$M^s_E := \{[\g] \in M^s : \E(\g) < E\}.$$
Geometrically, Theorem \ref{ThmSmallEnergy} says that for sufficiently small energy $\eps$, the space $M^s_\eps$ is what we expect it to be in light of the Morse-Bott analysis of the previous section.  Namely, $M^s_\eps$ is an open neighborhood of the the critical set of our flow, the space of vortices, within the stable manifold of the flow.  (Since we are working modulo gauge, the stable manifold in question is with respect to the flow on some coordinate patch near a vortex, as we analyzed in the previous section.)  Thus, while $M^s_\eps$ is an infinite-dimensional Hilbert manifold, the only topologically nontrivial portion of it comes from the finite dimensional space of vortices over which it fibers.  Furthermore, the Seiberg-Witten flow provides a weak homotopy equivalence from the entire space $M^s$, whose exact nature we do not know, to the small energy space $M^s_\eps$.

\begin{Theorem} \label{ThmSmallEnergy} (Small Energy Moduli Space)
We have the following:
 \begin{enumerate}
     \item For every $E > 0$, the inclusion $M_E^s \hookrightarrow M^s$ induces a weak homotopy equivalence.
     \item There exists an $\eps_0 > 0$ such that for all $0 < \eps < \eps_0$, the space $M^s_\eps$ is diffeomorphic to a Hilbert ball bundle over the vortex moduli space $\V(\Sigma)$.
 \end{enumerate}
\end{Theorem}

\Proof
(i) We want to show that the inclusion induces an isomorphism on all homotopy groups.  For this, we only have to show that $M_E^s \hookrightarrow M^s$ is surjective on all homotopy groups.  So let $f: S_n \to M^s$ be a representative element of $\pi_n(M^s)$ for some $n$.  Observe that for every $T \geq 0$, we have a continuous map $\tau_T:M^s \to M^s$ which translates an element by time $T$, i.e. $\tau_T(\g) = \g(\cdot + T)$.  Since the image of $f(S_n)$ is compact, and because energy is continuous on $M^s$, it follows that we can find a large $T$ such that $\tau_T(f(S_n)) \subset M^s_E$.  Thus, $\tau_{t}$, $0 \leq t \leq T$, provides a homotopy from $\tau_T(f(S_n))$ to $f(S_n)$.  Since $f: S^n \to M^s$ was arbitrary, this proves the desired surjectivity of the inclusion map on homotopy groups.

(ii) The Chern-Simons-Dirac functional $CSD^\Sigma$, being a Morse-Bott functional on the quotient space $\fB^{s-1/2}_{\eta^0}(\Sigma)$, is a small lower order perturbation of a positive-definite quadratic form when restricted to small neighborhood of the stable manifold to a critical point.  Hence, the level sets of energy on such a stable manifold, for energy close to the energy of the critical set, are just smooth spheres. Thus, the union of those level sets of energy less than $\eps$, which is precisely $M^s_\eps$, forms a Hilbert ball bundle over $\V$.  Here, in this last statement, we implicitly used Lemma \ref{LemmaSmallEnergy}, which tells us that for small enough energy, every configuration is gauge equivalent to a path that remains in a small $H^{1/2}(\Sigma)$ neighborhood of $\V$ for all time, in which case the above local analysis of $CSD^\Sigma$ near its critical set applies.\End

From the previous theorems, we can deduce the following theorem, which allows us to obtain Lagrangian submanifolds of $\fB^{s-1/2}_{\eta^0}(\Sigma)$ whose topology we can understand.  Namely, we consider the initial data of configurations in $\B^{s-1/2}_{\eta^0}(\Sigma)$ that converge under the Seiberg-Witten flow to a submanifold $\mathscr{L}$ inside the vortex moduli space $\V$ at infinity.  More precisely, define
the space
$$M^s_{\sL} = \{[\g] \in M^s : \partial_\infty[\g] \in \sL\}.$$
of monopoles in $M^s$ that converge to $\sL$.  For any $E>0$, we can also define
$$M^{s;\delta}_{\sL, E} = M^s_{\sL} \cap M^s_E.$$

\begin{Theorem}\label{ThmLagVort}
Let $\mathscr{L} \subset \V(\Sigma)$ denote any Lagrangian submanifold.
 \begin{enumerate}
   \item The space $M^s_\sL$ can be given the topology of a smooth Hilbert manifold.  The map $r_\Sigma: M^s_\sL \to \fB^{s-1/2}_{\eta^0}(\Sigma)$ is a diffeomorphism onto a Lagrangian submanifold of $\fB^{s-1/2}_{\eta^0}(\Sigma)$.  The space $M^s_\sL$ is weakly homotopy equivalent to a Hilbert ball bundle over $\sL$.
   \item If $\eps > 0$ is sufficiently small, then (i) holds with $M^s_{\sL,\eps}$ in place of $M^s_\sL$, and with ``weakly homotopy equivalent" replaced with ``diffeomorphic".
 \end{enumerate}
\end{Theorem}

\Proof Since the map $\partial_\infty: M^s \to \V$ is a smooth submersion, it follows $M^s_\sL \subset M^s$ has the topology of a smooth Hilbert manifold.  From Theorem \ref{ThmFinEnergy}(ii), we see that given $[\g] \in M^s_\sL$, the space $r_\Sigma(T_{[\g]}M^s_\sL)$ yields a Lagrangian subspace inside the symplectically reduced space
$$r_\Sigma(T_{[\g]}M^s)/r_\Sigma(T_{[\g]}M^s({\partial_\infty[\g]})).$$
This shows that $r_\Sigma(T_{[\g]}M^s_\sL)$ is a Lagrangian subspace of $T_{r_\Sigma[\g]}B_\mu^{s-1/2}(\Sigma)$. The remaining statements are now immediate.\End

Of course, having worked initially in the Hilbert space setting (as is necessary to use Banach space methods), one can then restrict to just those configurations that are smooth.  Thus, all the results above carry over mutatis mutandis to the smooth setting. In what follows, omission of the superscript $s$ from a configuration space denotes we are considering those configurations that are smooth.


\subsection{The General Case}

We now consider the general case of a $3$-manifold $Y$ with cylindrical ends.  Write $Y = Y_0 \cup ([0,\infty) \times \Sigma)$ as the union of a compact $3$-manifold with boundary $Y_0$ and a cylindrical end on which all structures (the metric, $\spinc$ structure, and perturbation) are product.  Using the previous results concerning monopole moduli spaces on semi-infinite cylinders together with those results of \cite{N1} on compact $3$-manifolds, we are able to describe the finite energy moduli space of monopoles on $Y$.  Indeed, on the compact part, we have the following theorem from \cite{N1}.  Given $\eta \in \Omega^1(Y_0; i\R)$ coclosed, we have the moduli space $M_\eta^s(Y)$ of $H^s(Y_0)$ gauge-equivalence classes of $\eta$-perturbed monopoles on $Y_0$.

\begin{Theorem}\label{ThmN1}
  Let $\s$ be a $\spinc$ structure on $Y_0$ such that $c_1(\s) \neq \frac{i}{\pi}[*\eta]$.  Then $M_\eta^s(Y_0)$ is a smooth Hilbert manifold and $r_\Sigma: M_\eta^s(Y_0) \to \fB_{\eta^0}^{s-1/2}(\Sigma)$ is a smooth submersion onto a Lagrangian submanifold.  The fiber of the submersion is isomorphic to the lattice $H^1(Y,\Sigma; \Z)$.
\end{Theorem}

Let $\eta \in \Omega^1(Y; i\R)$ be a coclosed one-form such that its restriction to the cylindrical end $[0,\infty) \times \Sigma$ is of the form $\eta^0 dt$ with $\eta^0$ time-independent.  Then $\eta$ induces for us perturbations $\eta|_{Y_0}$ and $\eta|_{[0,\infty)\times\Sigma}$ to the monopole equations on $Y_0$ and $[0,\infty) \times \Sigma$.  We thus obtain the corresponding perturbed moduli space of smooth monopoles:
\begin{align}
  M_\eta([0,\infty)\times\Sigma) &= M_{\eta|_{[0,\infty)\times\Sigma}}([0,\infty)\times\Sigma))\\
  M_\eta(Y_0) &= \{u \in \fC(Y_0):  SW_3(u) = \eta|_{Y_0}\}/\G(Y_0)\\
  M_\eta(Y) &= \{u \in \fC^\delta(Y) : SW_3(u) = \eta,\; \E(\g|_{[0,\infty) \times \Sigma_i}) < \infty,\; 1 \leq i \leq n\}/\G^\delta(Y).
\end{align}
From now on, we will not always distinguish between a form $\eta$ on $Y$ and its restriction to smaller domains in our notation as in the above.

As in \cite{KM}, we can describe $M_\eta(Y)$ as a fiber product of the moduli space of monopoles on $Y_0$ and on the ends.  We have restriction maps to the quotient configuration space $\fB(\Sigma)$ on the interface $\Sigma = \{0\} \times \Sigma$ of $Y_0$ and $[0,\infty) \times \Sigma$:
\begin{align}
  r_\Sigma^+: \fB(Y_0) & \to \fB(\Sigma)\\
  r_\Sigma^-: \fB([0,\infty)\times\Sigma) & \to \fB(\Sigma).
\end{align}
Via restriction, these maps then give us maps
\begin{align}
  r_\Sigma^+: M_\eta(Y_0) & \to \fB_{\eta^0}(\Sigma) \label{eq2:r+}\\
  r_\Sigma^-: M_\eta([0,\infty)\times\Sigma) & \to \fB_{\eta^0}(\Sigma). \label{eq2:r-}
\end{align}
One can show, as in \cite[Lemma 24.2.2]{KM}, the following:

\begin{Lemma}\label{LemmaFP}
The natural map $M_\eta(Y) \to M_\eta(Y_0) \times M_\eta([0,\infty)\times\Sigma)$ given by restriction yields a homeomorphism  from $M_\eta(Y)$ onto the fiber product of (\ref{eq2:r+}) and (\ref{eq2:r-}).
\end{Lemma}

Note that this lemma requires that there be no reducibles, which we always take to be the case for a suitable choice of perturbation and $\spinc$ structure landing us in the Morse-Bott situation.

Our main result is the following.  Pick a $\spinc$ structure $\s$ and perturbation $\eta$ on $Y$ as above and write $\Sigma = \partial Y_0$ as $\Sigma = \coprod_i \Sigma_i$ in terms of its connected components.  Let $d_i = \frac{1}{2}\left<c_1(\s),\Sigma_i\right>$ and where $g_i$ is the genus of $\Sigma_i$. Then the vortex moduli space at infinity can be written as
\begin{equation}
  \V(\Sigma) = \prod_{i=1}^n \V(\Sigma_i) \label{prodvort}
\end{equation}
where the precise degree of these vortex moduli spaces are given by Lemma \ref{LemmaVortices}.  Here, we suppose that $\eta$ is such that
\begin{equation}
  \frac{i}{\pi}\int_{\Sigma_i}*\eta \neq d_i \label{goodpert}
\end{equation}
for every $i$, so that our previous Morse-Bott analysis applies on each end. This is not always possible. Indeed, if $Y$ has a single end, both sides of (\ref{goodpert}) are always zero. However, (\ref{goodpert}) is possible  if $Y$ has at least two ends (here we assume $Y$, without loss of generality, is connected). For in this case, the restriction of $H^2(Y_0)$ to each boundary component of $Y_0$ has nontrivial image. Furthermore, we endow (\ref{prodvort}) with the product symplectic structure (of course, weighted with signs according as to whether the $\Sigma_i$ are incoming or outgoing ends, which we suppress from the above notation).

We can define the smooth map
\begin{align}
  \partial_\infty: M_\eta(Y) \to \prod_{i=1}^n \V(\Sigma_i) \label{eq2:partialM}
\end{align}
which sends a monopole to the gauge-equivalence class of its limit on each end.

\begin{Theorem}\label{ThmOne}
  Fix a coclosed perturbation $\eta \in \Omega^0(Y; i\R)$ as above and suppose it can be chosen so that (\ref{goodpert}) holds.  Then for a generic coclosed perturbation $\bar\eta$ compactly supported in the interior of the collar neighborhood $[-1,0] \times \Sigma \subset Y_0$, the space $M_{\eta+\bar\eta}(Y)$ is a smooth, compact, orientable manifold.   Moreover, the map (\ref{eq2:partialM}) is a Lagrangian immersion.
\end{Theorem}

To prove the first part of the theorem, we show that for generic $\bar\eta$, the restriction maps given by (\ref{eq2:r+}) and (\ref{eq2:r-}) intersect transversally.  This requires understanding how the parametrized moduli spaces $M_{\eta+\bar\eta}(Y)$, as $\bar\eta$ varies over a Banach space of coclosed forms with compact support in the interior of $Y_0$, restrict to the boundary $\Sigma$.

In detail, fix a small open interval $I \subset [-1,0]$ whose closure is contained within $[0,1]$.  Fix a larger interval $\tilde I \supset I$ with the same property.  Fix a countable collection of smooth compactly coclosed imaginary $1$-forms supported in $\tilde I \times \Sigma$ such that their restrictions to $I \times \Sigma$ are dense in the space of smooth coclosed $1$-forms on $I \times \Sigma$.  Then as in \cite{KM}, one can form a Banach space of smooth\footnote{One could work with a Banach space of $H^s$ forms for large $s$, but for convenience, we will take our Banach space to consist of smooth forms.} coclosed forms $\P$ that is given by the closure of the span of the given countable collection.  We will take $\P$ to be our space of perturbations to the Seiberg-Witten equations, where an element of $\P$ is said to be generic if, as usual, it lies within some unspecified residual subset of $\cP$.

\begin{Lemma}\label{LemmaPert}
  Let $s \geq 2$. Consider the map
  \begin{align*}
    \F: \fC^s(Y_0) & \times \cP \to \K^{s-1}(Y_0) \times \fB^{s-1/2}(\Sigma)\\
    (u,\bar\eta) & \mapsto (SW_3(u) - (\eta + \bar\eta), r_\Sigma(u)),
  \end{align*}
  where $\K^{s-1}(Y_0)$ denotes the Hilbert bundle over $\fC^s(Y_0)$, whose fiber over $u$ is $\K_u^{s-1}$.  Then for any $(u,\bar\eta)$ such that $SW_3(u) = \eta + \bar\eta$, the image of $D_{u,\bar\eta}F$, orthogonally projected into $\K_u \oplus T_{r_\Sigma(u)}\mu^{-1}(\eta^0)$, is dense with respect to the $L^2$ topology.
\end{Lemma}

\begin{proof}
  Let $(v_1,v_2) \in \K^{s-1}_u \times \T_\Sigma^{s-1/2}$ be $L^2$ orthogonal to the image of $D_{u,\bar\eta}F$.  This means
  $$(\H_u(\delta u) + \delta\bar\eta, v_1)_{L^2(Y_0)} + (r_\Sigma(\delta u), v_2)_{L^2(\Sigma)} = 0$$
  for all $(\delta u,\delta\bar\eta) \in T_u\fC^s(Y_0) \times \cP$.  Let $\delta\bar\eta = 0$.  Choosing $\delta u$ compactly supported, integration by parts shows that $v_1 \in \ker \H_u$.  Then choosing $\delta u$ arbitrary, integration by parts again allows us to conclude that $v_2 = Jr_\Sigma(v_1)$.  Letting $\delta u = 0$ and let $\delta\bar\eta$ vary.  Then on $U := I \times \Sigma$, writing ${v_1}|_U = (b, \Phi)$, where $b \in \Omega^1(U; i\R)$ and $\Phi$ is a spinor, it follows that $b$ is closed.  On the other hand, since $u =: (B,\Psi)$ is an irreducible monopole, then $\Psi$ vanishes only on a set of isolated points by unique continuation.  Thus, since $(d\xi, \Phi) \in \ker \H_\g$ implies $\mathrm{Im}\rho^{-1}(\Phi\Psi^*)_0 = 0$, we must have $\Phi = \xi\Psi$ for some imaginary valued function $\xi$.  Then, $(b, \xi\Psi) \in \ker \H_u$ implies $\rho(b)\Psi + \rho(d\xi)\Psi = 0$, since $D_B\Psi = 0$.  Thus, we see that ${v_1}|_U$ is of the form $(-d\xi, \xi\Psi)$, i.e. it lies in the infinitesimal gauge orbit through $(B,\Psi)$.  Placing $v_1$ in temporal gauge along the entire collar neighborhood $[-1,0] \times \Sigma$, unique continuation (see \cite{KM} or \cite{N1}) implies ${v_1}|_{[-1,0] \times \Sigma}$ equals $(-d\xi,\xi\Psi)$ for some $\xi \in \Omega^0([-1,0] \times \Sigma)$.  It follows that ${v_2} = Jr_\Sigma {v_1}$ belongs to $J\J_{r_\Sigma(u)}$.  But $J\J_{r_\Sigma(u)} \cap T_{r_\Sigma (u)}\mu^{-1}(\eta^0) = 0$ by Proposition \ref{PropMomentMap}, so $v_2 = 0$.  Thus $r_\Sigma v_1 = 0$ and again by unique continuation $v_1 = 0$.  This proves that the $L^2$ closure of the projection of the range of $D_\g F$, projected onto $\K_u \oplus T_{r_\Sigma(u)}\mu^{-1}(\eta^0)$, is equal to $\K_u \oplus T_{r_\Sigma(u)}\mu^{-1}(\eta^0)$.\End
\end{proof}

\begin{Corollary}\label{Cor5.11}
  Let $\cZ = \{([u],\bar\eta) \in \fB^s(Y_0) \times \cP : SW_3(u) = \eta+\bar\eta\}$.  Then $\cZ$ is a smooth Banach manifold and the maps $r_\Sigma^-: \cZ \to \B_{\eta^0}^{s-1/2}(\Sigma_0)$ and $r_\Sigma^+: \M_\eta([0,\infty)\times\Sigma) \to \B_{\eta^0}^{s-1/2}(\Sigma_0)$ are transverse.
\end{Corollary}

\Proof For any $\bar\eta$, at any point of the fiber product of the maps $r_\Sigma^-: \M_{\eta + \bar\eta}^s(Y_0) \to \B_{\eta^0}^{s-1/2}(\Sigma)$ and $r_\Sigma^+: M_\eta^s([0,\infty)\times\Sigma) \to  \B_{\eta^0}^{s-1/2}(\Sigma)$, the image of the differentials of these maps yield a Fredholm pair of Lagrangian subspaces.  Moreover, these subspaces remain Fredholm if we pass to the $L^2$ closure.  Next, it is an elementary fact that if a subspace of a Banach space has finite codimension and is dense then it must be the entire space.  The previous lemma now implies the corollary.\End

\textbf{Proof of Theorem \ref{ThmOne}:} By Corollary \ref{Cor5.11}, the fiber product of the maps $r_\Sigma^-: \cZ \to \fB^{s-1/2}_{\eta^0}(\Sigma_0)$ and $r_\Sigma^+: \M_\eta([0,\infty)\times\Sigma) \to \fB^{s-1/2}_{\eta^0}(\Sigma_0)$ is a smooth submanifold of $\cZ \times \M_\eta([0,\infty)\times\Sigma)$.  Moreover, because $r_\Sigma^+(M_{\eta+\bar\eta}(Y_0))$ and $r_\Sigma^-(M_{\eta}([0,\infty)\times\Sigma))$ always intersect in a Fredholm manner, the projection of the above fiber product onto the space of perturbations $\cP$ is Fredholm.  By the Sard-Smale theorem, we can find a residual set of regular values for this projection.  Choosing $\bar\eta$ to be such a regular value, we then obtain that the restriction maps from the corresponding moduli spaces are transverse.

It follows that the resulting fiber product $M_{\eta+\bar\eta}(Y)$ is smooth and finite-dimensional.  The fact that $M_{\eta+\bar\eta}(Y)$ is compact follows from the compactness results for the perturbed Seiberg-Witten equations, see \cite[Chapter 24.5]{KM}.  In our situation, all finite energy configurations in $M_{\eta+\bar\eta}(Y)$ must have the exactly the same (perturbed) topological energy, since the space of vortices on each end is connected (and so $CSD^{\Sigma_i}$ has constant value on the vortices on each $\Sigma_i$).  Moreover, we cannot have trajectory breaking on the ends for the same reason: the only finite energy solutions on an infinite cylinder $(-\infty,\infty) \times \Sigma_i$ are translation-invariant zero energy vortices.  Thus, our space $M_{\eta+\bar\eta}(Y)$ is compact as is.

For the second statement, we can see this very easily in geometric terms.  From Theorem \ref{ThmN1}, we know that the image of (\ref{eq2:r+}) is a Lagrangian submanifold. Let $[\g] \in M_\eta([0,\infty)\times\Sigma)$ and define $[\a] := \partial_\infty[\g] \in \V(\Sigma)$.  Note that the differential of $\partial_\infty$ at $[\g] \in M_\eta([0,\infty)\times\Sigma)$ has kernel precisely equal to $T_{[\g]}M_{\eta}([\a])$, the tangent space to the stable manifold to $[\a]$.  On the other hand, by Theorem \ref{ThmFinEnergy}(ii), we have that $r_\Sigma^- \big(T_{[\g]}M_\eta([\a]))$ is an isotropic subspace annihilating the coisotropic subspace $r_\Sigma^- \big(T_{[\g]}M_\eta([0,\infty)\times\Sigma)\big)$.  So given any $u \in M_{\eta+\bar\eta}(Y)$, it follows that the differential
$$D_{u}\partial_\infty: T_{u}M_{\eta+\bar\eta}(Y) \to T_{\partial_\infty (u)}\V(\Sigma)$$
has range isomorphic to the symplectic reduction of the Lagrangian subspace $T_{r_\Sigma(u)}\big(r_\Sigma^+ M_{\eta+\bar\eta}(Y_0)\big)$ coming from $Y_0$ with respect to the coisotropic space $T_{r_\Sigma(u)}r_\Sigma^- (M_\eta([0,\infty)\times\Sigma)\big)$ coming from the ends.  In particular, the differential of $\partial_\infty$ at any monopole on $Y$ has image a Lagrangian subspace.  Moreover, the map $\partial_\infty$ is an immersion due to the transversality of the maps (\ref{eq2:r+}) and (\ref{eq2:r-}), which implies that the symplectic reduction in question is injective. This proves the theorem.\End

\begin{Remark}
  Our proof consisted of patching together analysis from the compact piece $Y_0$ with the cylindrical end $[0,\infty) \times \Sigma$.  One could work directly on $Y$, proving that the map $SW_3: \fC^\delta(Y) \to \K^\delta(Y)$ has $\eta + \bar\eta$ as a regular value for generic $\bar\eta$.  Since the linearization of $SW_3$ is formally self-adjoint, we know that the image of $(\partial_\infty)_*$ on tangent spaces of $\M_{\eta+\bar\eta}(Y)$ are isotropic, and by the Atiyah-Patodi-Singer index theorem, we know that $\dim \M_{\eta+\bar\eta}(Y) = \frac{1}{2}\dim \V(\Sigma)$.  However, we would not know a priori that $\partial_\infty: \M(Y) \to \V(\Sigma)$ were an immersion.  Studying how one can choose $\bar\eta$ so that $\M_{\eta+\bar\eta}(Y)$ is an immersed moduli space seems like it would involve the same amount of work as carried out above.  What makes the above approach, involving the fiber product description of $\M(Y)$, somewhat miraculous is that one gets the immersive property automatically from the transverse intersection occurring in a symplectic reduction.  That is, the immersive property is obtained for free from the symplectic geometry.
\end{Remark}

Next, we show that the moduli space $M(Y)$ is an orientable manifold.  In fact, we show that the determinant line bundle $\Lambda(Y)$ over $\B^\delta(Y)$, whose fiber over $[u] \in \B^\delta(Y)$ is $\det(\bH_u)$, is trivial.  Thus, an orientation for $\Lambda$ determines a unique orientation of $M(Y)$, which a priori may have many since $M(Y)$ could be disconnected.

\begin{Lemma}\label{LemmaOr}
  The determinant line bundle $\Lambda(Y)$ is trivial.
\end{Lemma}

\Proof We know that $\B^{s,\delta}(Y)$ is diffeomorphic to a Hilbert bundle over $\V_Y(\Sigma)$, and hence homotopy equivalent to this latter space, with the homotopy equivalence being given by the limiting boundary value map $\B^\delta(Y) \to \V_Y(\Sigma)$.  Thus, suppose we are given an arbitrary loop $z_0: S^1 \to \V_Y(\Sigma)$.  We will extend it to a loop $z: S^1 \to \B^\delta(Y)$ and show that $\Lambda(Y)$ restricts trivially to $z$. Since $z_0$ is arbitrary, this will show that $\Lambda(Y)$ is a trivial line bundle, since it restricts trivially to every homotopy class of loops in $\B^{s,\delta}(Y)$. So given $z_0$, we construct an extension $z$ piece by piece as follows.  Extend $z_0$ along the neck $[1,\infty)\times\Sigma$ in a time-independent manner.  The homotopy $(C,\Psi) \to (C,t\Psi)$, $t \in [0,1]$ is a gauge-equivariant homotopy inside $\fC(\Sigma)$, and we can apply this homotopy to obtain a path of configurations $z: S^1 \to \fB^\delta([0,\infty) \times \Sigma)$, with each configuration extended to $[0,1]\times\Sigma$ via the above homotopy.  We can then extend $z$ in some smooth manner into the interior of $Y$, with $z$ reducible on $Y_0$, thus obtaining our path $z: S^1 \to \fB^\delta(Y)$. 

Orienting the determinant line $\Lambda(Y)$ on this loop $z$ is not straightforward since we do not have a complex linear family of operators on $Y$.  We want to pass to a situation in which this is true however.  To do this, we use an excision argument to pass to operators defined on cylinders where complex linearity can be exploited.  Consider the loop $z'(\theta)$, $\theta \in S^1$, of configurations on $Y$ such that $z'(\theta)|_{Y_0} = z(\theta)|_{Y_0}$ and $z(\theta)|_{[0,\infty)\times\Sigma}$ is the time-translation invariant reducible configuration $z(\theta)|_{\{0\} \times \Sigma}$.  In this way, we obtain a loop $z'$ of reducible configurations on $Y$ that agrees with $z$ on $Y_0$.

On the other hand, we can also define loops of configurations $\tilde z$ and $\tilde z'$ on $\R \times \Sigma$, which agree with $z$ and $z'$ on $[0,\infty)\times\Sigma$ and which are time-translation invariant extensions of their values on $\{0\}\times\Sigma$ to $(-\infty,0]\times\Sigma$.  By standard excision properties of determinant lines \cite[Chapter 7]{DK}, \cite[Chapter 20.3]{KM}, we have
\begin{equation}
  \det (\bH_z) \otimes \left(\det (\bH_{z'}) \right)^* \cong \det(\bH_{\tilde z}) \otimes \left(\det (\bH_{\tilde z'})\right)^* \label{exc}
\end{equation}
where the left-hand side consists of a determinant line of a loop of operators over $Y$ and the right-hand side consists of a determinant line of a loop of operators over $\R \times \Sigma$. Here, care must be taken in our choice of initial loop $z_0: S^1 \to \V_Y(\Sigma)$ because this loop determines the asymptotic behavior of the loops $z'$, $\tilde z$, and $\tilde z'$. We must ensure that each corresponding loop of Hessian operators occuring in (\ref{exc}) parametrizes a family of Fredholm operators, which means that we must be able to choose a suitable weighted space for the domain and range of the loop of operators which renders them all Fredholm.

We show that for a suitable representative of $z_0$ and negative weight $-\delta$ on the cylindrical ends of $Y$ and of $\R \times \Sigma$, with $\delta > 0$ sufficiently small, all the operators appearing in (\ref{exc}) define a family of Fredholm operators. For $\bH_z$ there is nothing to show, since the Morse-Bott assumption tells us precisely that for sufficiently small negative weight $\delta$, no eigenvalues of the linearized gauge-fixed vortex operator $\barB_{z_0(t)}$ can cross $-\delta < 0$ for $\delta$ a sufficiently small weight. To prove that $\bH_{z'}$ defines a family of Fredholm operators, we need to prove that no eigenvalues for the family of operators $\barB_{z'_0(t)}$ cross $-\delta$, where $z'_0(t)$ is a loop of connections obtained from $z_0(t)$ by setting the spinor equal to zero.

In this case, the family $\barB_{z'_0(t)}$ is the direct sum of a single Hodge operator and a family of Dirac operators parametrized by $z'_0(t)$. We want this family of Dirac operators to have kernel of constant dimension. This will be true as long as the loop $z'_0(t)$, viewed as a loop in the Jacobian variety of $\Sigma$, avoids the locus of line bundles that are determined by special divisors. More precisely, we wish for the operators
\begin{align*}
  \partial_{z'_0(t)}: (K \otimes L)^{1/2} & \to (K^{-1} \otimes L)^{1/2}\\
  \bar\partial_{z'_0(t)}^*: (K^{-1} \otimes L)^{1/2} & \to (K \otimes L)^{1/2}
\end{align*}
to have kernel of constant dimension for all $t \in S^1$. This is guaranteed if the family of divisors corresponding to the holomorphic structures induced on $(K \otimes L)^{1/2}$ and $-(K^{-1} \otimes L)^{1/2}$ by $z'_0(t)$ avoid those which are special. Thus, the loop $z_0'$ and hence $z_0$ must avoid a proper analytic subvariety, which is at least a complex codimension one condition. Since a loop is of one real dimension, for any given homotopy class of loops we can choose a representative $z_0$ which avoids the subvariety.

The exact same analysis applies at the negative infinite end of $\R \times \Sigma$. This shows that we can choose a small negative weight for which the operators in (\ref{exc}) are all Fredholm, which justifies the excision procedure.

From (\ref{exc}), to orient $\det (\bH_z)$, it suffices to orient the other three lines.  The lines $\det (\bH_{z'})$ and $\det (\bH_{\tilde z'})$ are trivial since they are families of reducible configurations (so that one has a fixed Hodge operator and family of complex linear Dirac operators).  It remains to orient the line $\det(\bH_{\tilde z})$. Here, we invoke Lemma \ref{LemmaZero}(iv), which tells us that the operators belonging to the family $\bH_{\tilde z}$ are all complex-linear with respect to $I$.  Thus, $\det (\bH_{\tilde z})$ is canonically oriented. Altogether, this shows that $\det (\bH_z)$, the line bundle $\Lambda(Y)$ restricted to the loop $z$, is trivial.  \End

\begin{Corollary}\label{CorOr}
  The moduli space $M(Y)$ is orientable with a unique orientation induced from an orientation for $\Lambda(Y)$. An orientation for $\Lambda(Y)$ is determined by an orientation for the image of $H^1(Y) \to H^1(\Sigma)$.
\end{Corollary}

\begin{Proof}
  We need only prove the latter statement. From (\ref{exc}), to orient $\Lambda(Y)$, it suffices to orient a single line $\left(\det (\bH_{z'(t_0)}) \right) \otimes \det(\bH_{\tilde z(t_0)}) \otimes \left(\det (\bH_{\tilde z'(t_0)})\right)^*$ for some $t_0$. The middle term arises from a complex linear operator and so has a canonical orientation. The first and last operators above are each given by reducible configurations $z'(t_0)$ and $\tilde z'(t_0)$. For a reducible configuration, we need only orient the determinant line of the associated Hodge operator, since the associated Dirac operator is complex linear. It follows that we need to orient the determinant line of the Hodge operator $(*d + d^*) \oplus d$ on $Y$ and $\R \times \Sigma$ with a small negative weight on each end. In general, on an arbitrary $3$-manifold $X$ with boundary, if we form $X^*$ by attaching cylindrical ends to $\partial X$ and choose a small negative weight on each end, the kernel of the Hodge operator is isomorphic to $H^1(X) \oplus H^0(X)$ (with real coefficients) while the cokernel is isomorphic to the image of $H^1(X,\partial X) \to H^1(X)$ (\cite{APS}, \cite{Don}). Thus, via exactness of $H^1(X,\partial X) \to H^1(X) \to H^1(\partial X)$ and since $H^0(X)$ has a canonical orientation, orienting the determinant line of the Hodge operator on $X^*$ with a small negative weight is equivalent to orienting the image of $H^1(X) \to H^1(\partial X)$. We now apply this result to $X = Y$ and $X = \R \times \Sigma$, and note that the image of $H^1(\R \times \Sigma) \to H^1(-\Sigma \times \Sigma)$ is isomorphic to $H^1(\Sigma)$, which has a canonical orientation since it is a symplectic vector space.\End
\end{Proof}


\section{Donaldson's ``TQFT"}

Consider a closed $3$-manifold $Y$ with $b_1(Y) > 0$.  Then it is possible to choose a connected nonseparating orientable hypersurface $\Sigma \subset Y$.  We can then form the cobordism $W: \Sigma \to \Sigma$ obtained by removing $\Sigma$ from $Y$ and then form the corresponding cylindrical end manifold
$$W^* = W\cup \left((-\infty,0] \times \Sigma\right) \cup \left([0,\infty) \times \Sigma\right).$$
Here of course, we assume $W$ has a metric which is product in a neighborhood of the boundary so that $W^*$ is a smooth, Riemannian manifold.  The original manifold $Y$ is obtained from $W$ by identifying the two boundary components by a diffeomorphism $h: \Sigma \to \Sigma$.

Let $\s_0$ be a $\spinc$ structure on $W^*$ and $\eta \in \Omega^1(W^*; i\R)$ a coclosed form satisfying the usual product structure assumptions on the ends as in the previous sections.  We can identify $\s_0$ with its restriction (again  denoted $\s_0$) to $W$, and let $\Spinc(Y,\s_0)$ denote the set of all $\spinc$ structures on $Y$ obtained by closing up $W$ by the diffeomorphism $h$ and by using all possible gluing parameters (i.e. the inequivalent ways of identifying the $\spinc$ structures on $\Sigma$) as indexed by $\Gamma = H^1(\Sigma; \Z)/\left( H^1(W; \Z)|_\Sigma\right)$.  Following the program set out by Donaldson in \cite{Don99}, we can compute the Seiberg-Witten invariants of the closed manifold $Y$, or more precisely, those obtained by summing over the $\spinc$ structures of $\Spinc(Y,\s_0)$ for some $\s_0$, using a topological quantum field-theoretic framework. We will refer to such a framework as a TQFT for brevity, even though it is not so on the nose as we shall see.

Let us first state the invariant this TQFT produces. As a set, the monopoles on $Y$ are precisely those monopoles on $W$ whose boundary values on the two components agree when we glue by the diffeomorphism $h$ and any gluing parameter.  However, to obtain an intersection problem involving finite dimensional objects (the space of monopoles on a compact $3$-manifold with boundary and the space of its boundary values are infinite dimensional modulo gauge), we instead attach cylindrical ends and consider the space of perturbed monopoles $M_\eta(W^*, \s_0)$ with respect to the $\spinc$ structure $\s_0$. Based on gluing principles, one expects that the moduli space of monopoles on $Y$, with respect to the $\spinc$ structures in $\Spinc(Y,\s_0)$ should correspond bijectively to the intersection of $\partial_\infty M_\eta(W^*, \s_0)$ with $\Gamma_h$, the graph of $h$, inside $\V(-\Sigma) \times \V(\Sigma)$.  Here a sum of $\spinc$ structures is involved since the vortex moduli space on $\Sigma$ is formed by dividing by gauge transformations on $\Sigma$, not all of which extend to $Y$.  Those that do not extend contribute to a gluing parameter that changes the glued $\spinc$ structure we obtain on $Y$, which results in a summation over elements of $\Spinc(Y,\s_0)$. Finally, to compute the Seiberg-Witten invariants of $Y$, one also needs count the monopoles on $Y$ with the appropriate signs.

The main result of this paper, confirming the picture outlined by Donaldson, is that this signed count corresponds precisely to the (homological) signed intersection of $\partial_\infty M(W^*, \s_0)$ with $\Gamma_h$.  The precise result is as follows.  Let $\eta \in \Omega^1(Y; i\R)$ be a coclosed form on $Y$ serving as a perturbation for the Seiberg-Witten equations on $Y$.  The effect that the perturbation has on the Seiberg-Witten moduli space only depends on the cohomology class of $*\eta$. By standard Hodge theory, one can always choose $\eta$, while remaining in a fixed cohomology class $[*\eta]$, so that in a tubular neighborhood $[-1,1]\times\Sigma$ of $\Sigma$, $\eta$ is of the form $\eta^0 dt$, with $\eta^0 \in \Omega^0(\Sigma; i\R)$ independent of the normal coordinate $t$.
Supposing $\eta$ is of that form then, it has a natural extension, in a time-translation invariant fashion, to the cylindrical end manifold $W^*$, which we again denote by $\eta$.

We say that a coclosed $\eta \in \Omega^1(Y; i\R)$ as above is \textit{admissible}.  Observe we can always chose an admissible $\eta$ so that $\left<*\eta, \Sigma\right>$ takes on any desired value, since $[\Sigma]$ is nontrivial.  Thus, from Lemma \ref{LemmaVortices}, we can always choose an admissible perturbation that places us within a Morse-Bott context.

\begin{Theorem}\label{ThmTwo}
  Let $Y$ be a closed $3$-manifold with $b_1(Y) > 0$.  Pick any connected nonseparating orientable hypersurface $\Sigma \subset Y$ and form the cylindrical end manifold $W^*$ from the manifold $W = Y\setminus\Sigma$ as above. Pick a $\spinc$ structure $\s_0$ on $W$ and fix an admissible perturbation $\eta \in \Omega^1(Y; i\R)$ such that $\frac{i}{\pi}\left<*\eta,\Sigma\right> \neq \frac{1}{2}\left<c_1(\s_0),\Sigma\right>$.  Then for generic coclosed perturbations $\bar\eta$ compactly supported in $W$, we have
  \begin{equation}
    \sum_{\s \in \Spinc(Y,\s_0)} SW(\s, \eta) = (\partial_\infty)_*[M_{\eta + \bar\eta}(W^*, \s_0)] \cap [\Gamma_h]. \label{intersect}
  \end{equation}
  Here, a homology orientation on $Y$ and an orientation of $M_{\eta + \bar\eta}(W^*, \s_0)$ are chosen compatibly (each of these determines an overall sign for the left-hand side and right-hand side, respectively).
\end{Theorem}

Since the Seiberg-Witten invariant vanishes for all but finitely many $\spinc$ structures, the above sum is well-defined.  Moreover, the sum only depends on the perturbation $\eta$ when $b_1(Y) = 1$.

\begin{Remark}\label{RemMark}
  Suppose $H^2(W) = \Z$.  Then the $\spinc$ structures occurring in (\ref{intersect}) are precisely those $\spinc$ structures $\s$ on $Y$ such that $\left<c_1(\s),\Sigma\right> = \left<c_1(\s_0),\Sigma\right>$.  This is the version that is implicitly being used in \cite{Mark}.
\end{Remark}

\begin{Remark}
  Let $d = \frac{1}{2}\left<c_1(\s),\Sigma\right>$ and $g = \mathrm{genus}(\Sigma)$. If $|d| > g-1$, the right-hand side of (\ref{intersect}) vanishes since we can choose $\eta$ so that the resulting vortex moduli spaces involved are empty by Lemma \ref{LemmaVortices}.  
\end{Remark}

As a simple application, we can recover the following well-known formulas for the Seiberg-Witten invariants of a product $3$-manifold.  In the case of $b_1(Y) = 1$, it is perhaps worth noting that no wall-crossing analysis is needed in our computation.

\begin{Corollary}
  Let $Y = S^1 \times \Sigma$.  Let $\s_d$ be the product $\spinc$ structure on $\Sigma$ such that $\frac{1}{2}\left<c_1(\s),\Sigma\right> = d$.
  \begin{enumerate}
    \item If $g \geq 1$, then
  $$SW(Y,\s_d) = \chi(\V_{g - 1 - |d|}(\Sigma)).$$
  \item Suppose $g = 0$.  Define
  $SW(Y,\s_d, \pm) = SW(Y,\s_d, \eta_\pm)$, where $\eta_\pm$ is any admissible perturbation such that $\pm\left(\frac{i}{\pi}\left<*\eta_\pm,\Sigma\right> - d\right) > 0$.  Then
  $SW(Y, \s_d, +) = d$ if $d \geq 0$ and zero otherwise, and $SW(Y,\s_d,-) = |d|$ if $d \leq 0$ and zero otherwise.
  \end{enumerate}
\end{Corollary}

\Proof This follows from (\ref{intersect}), the fact that the Seiberg-Witten invariants of $Y$ are only supported on product $\spinc$ structures, and Lemma \ref{LemmaVortices}. Here, $h = \mathrm{id}$ and we perturb so that the transverse intersection (\ref{intersect}) is simply the homological self-intersection of $\V(\Sigma)$ inside the diagonal $\V(\Sigma) \times \V(\Sigma)$.\End

\begin{Remark}\label{RemGen}
  Our methods actually establish a slightly more general formula than (\ref{intersect}).  Indeed, one can instead study Lagrangian intersections in some covering of the vortex moduli spaces on $\Sigma$ obtained by only dividing by some subgroup of the gauge group on $\Sigma$.  We would then have a sum over only a corresponding smaller set of gluing parameters on the left-hand side of (\ref{intersect}). However, this resulting sum of Seiberg-Witten invariants would not be obtained from Donaldson's TQFT picture we describe.
\end{Remark}

The above theorem has a natural TQFT-like formulation following \cite{Don99}.  This is because the above signed intersection can be regarded as a graded trace of a ``push-pull" map on the total homology of vortex moduli spaces, analogous to the trace formula that appears in the Lefschetz number of a self-map of a space.  Moreover, gluing allows us to decompose such a push-pull map into a composite of push-pull maps between intermediate vortex moduli spaces when we write our cobordism as a composite cobordism (satisfying a transversality hypothesis, see (\ref{Lemma-transglue})) and stretch along the intermediary neck joining the cobordisms.

More precisely, suppose we are given an arbitrary cobordism $W: \Sigma_0 \to \Sigma_1$, where as usual, we assume our cobordism to carry product structures (metric, $\spinc$, and perturbing coclosed $1$-forms of the type described above) near the boundary as needed.  Attaching cylindrical ends as before to obtain the cylindrical end manifold $W^*$, we can consider the moduli space $M_\eta(W^*)$ of $\eta$-perturbed finite energy monopoles on $W^*$. By Theorem \ref{ThmOne}, we know that
$$\partial_\infty: M_\eta(W^*) \to \V(-\Sigma_0) \times \V(\Sigma_1)$$
is a Lagrangian immersion.  Since $M_\eta(W^*)$ is a orientable by Corollary \ref{CorOr}, it carries a fundamental class which we may push-forward (in homology with real coefficients):
\begin{align}
  \rho_W := (\partial_\infty)_*[M_\eta(W^*)] & \in H_*(\V(-\Sigma_0) \times \V(\Sigma_1)) \nonumber \\
  & \cong H_*(\V(-\Sigma_0)) \otimes H_*(\V(\Sigma_1)) \nonumber \\
  & \cong H^*(\V(-\Sigma_0)) \otimes H_*(\V(\Sigma_1)) \nonumber \\
  & \cong \mathrm{Hom}(H_*(\V(-\Sigma_0)), H_*(\V(\Sigma_1)), \label{rhoW}
\end{align}
Here, we used the K\"unneth formula and Poincar\'e duality in the above. 

We thus have the following ``TQFT construction" of the Seiberg-Witten invariants of a closed oriented $3$-manifold $Y$ with $b_1(Y) > 0$. Fix two parameters $d \in \Z$ and $\underline{\eta} \in \R \setminus \Z$. To each Riemann surface $\Sigma$, we assign the graded vector space $H_*(\V(\Sigma))$ where $\V(\Sigma) = \V_k(\Sigma)$ is the degree $k$ vortex moduli space on $\Sigma$, where $k = k(\Sigma,d, \underline{\eta})$ as given by (\ref{kpm}) and Lemma \ref{LemmaVortices}, with $\underline{\eta}$ representing the value of $\frac{i}{2\pi}\int \eta^0$ in that lemma. To each \textit{elementary} cobordism $W: \Sigma_0 \to \Sigma_1$, we assign the morphism $\rho_W$, where the $\spinc$ structure $\s$ on $W$ is the one determined by requiring $\frac{1}{2}\left<c_1(\s),\Sigma_0\right> = d$.  Functoriality with respect to composition of transverse elementary cobordisms (see the next section) follows from the appropriate gluing results, which we will analyze soon.  The number we associate to a closed manifold $Y$ which is obtained by closing up the composite of transverse elementary cobordisms $W = W_n \circ \cdots \circ W_1: \Sigma \to \Sigma$ by a diffeomorphism $h$ is simply the graded trace of $h_* \circ \rho_{W_n} \circ \cdots \circ \rho_{W_1}: H_*(\V(\Sigma)) \to H_*(\V(\Sigma))$ (where the grading is the natural one on homology).  In fact, working through the definitions, this graded trace corresponds precisely to the homological intersection (\ref{intersect}).

Note that this construction is not quite a TQFT in several senses.  First, the invariant it computes for a closed manifold $Y$ is only a topological invariant (i.e. only depends on $d$) when $b_1(Y) > 1$.  Moreover, in this case, we only a priori know that this quantity is a topological invariant from the topological invariance of Seiberg-Witten theory itself.  Indeed, in the TQFT construction, we had to choose various compatible metrics and perturbations and there is no obvious reason, based on the TQFT definition alone, why the associated numerical invariant we obtain on a closed manifold should be independent of those choices (and indeed it is not in case $b_1(Y) = 1$).  If $b_1(Y) = 1$, choosing $d$ and $\underline{\eta}$ so that we land in the first case of Lemma \ref{LemmaVortices}, i.e. letting $k = k_+ + d$, means that the Seiberg-Witten invariants we compute is with respect to the chamber of $H^1(Y)$ determined by the ray $\lambda \mathrm{PD}([\Sigma])$, $\lambda > d$, where $\mathrm{PD}([\Sigma])$ denotes the Poincar\'e dual of $\Sigma$ . (Had we defined $k = k_- + d$ we would get the complementary chamber.)


Second, while the morphism $\rho_W$ is defined above for any cobordism, it depends on the choice of a $\spinc$ structure on $W$.  For an elementary cobordism $W$, $H^2(W; \Z) = \Z$ and thus the $\spinc$ structures are uniquely parametrized by the evaluation of their first Chern class along one of the boundary components.  Hence, if we wish to get a well-defined morphism that depends only on the fixed parameter $d$, we must work with elementary cobordisms.  Moreover, composability of our cobordisms requires a transversality hypothesis, see Definition \ref{DefTrans}.

Nevertheless, this TQFT like construction is a powerful point of view because the composition rule allows one to reduce the computation of the Seiberg-Witten invariants of a closed manifold to understanding how these push-pull maps behave on just elementary cobordisms. Recall that any cobordism, in particular, the one obtained from $Y\setminus\Sigma$, can be decomposed into a composite of elementary cobordisms.  Moreover, one can arrange this decomposition so that all the cobordisms are transverse, see \cite{Mark}.  Donaldson explicitly computes what the map $\rho_W$ is for an elementary cobordism using only elegant topological and algebraic arguments in \cite{Don99}.
For completeness, we describe these maps explicitly, following \cite{Mark}. Recall that as graded vector spaces, we have the isomorphism
\begin{equation}
  H^*(\Sym^k(\Sigma)) \cong \bigoplus_{i=0}^k \Lambda^i(H^1(\Sigma)) \oplus \mathrm{Sym}^{k-i}(H^0(\Sigma) \oplus H^2(\Sigma)) \label{eq:symprod}
\end{equation}
where the right-hand-side is graded in the natural way. Here, $\Lambda^i$ and $\mathrm{Sym}^i$ denote the $i$th exterior and symmetric powers, respectively.  An elementary cobordism $W: \Sigma_g \to \Sigma_{g+1}$ is given by attaching a $1$-handle to $[0,1] \times \Sigma_g$ at $\{1\}\times\Sigma$.  Let $c \in H^1(\Sigma_{g+1})$ be the cocycle Poincar\'e dual to the intersection of the cocore of the $1$-handle with $\Sigma_{g+1}$.  The cobordism $W$ gives us an imbedding $H^1(\Sigma_g) \hookrightarrow H^1(\Sigma_{g+1})$ and hence an imbedding of $H^*(\Sym^k(\Sigma_g)) \hookrightarrow H^*(\Sym^k(\Sigma_{g+1}))$ using (\ref{eq:symprod}).  Under this identification, the map $\rho_W$ is explicitly given by
\begin{align}
  \rho_W: H^*(\Sym^k(\Sigma_g)) & \to H^*(\Sym^{k+1}(\Sigma_{g+1}))\nonumber \\
  \omega & \mapsto c \wedge \omega.
\end{align}
If we reverse the cobordism and consider a $2$-handle attachment $W: \Sigma_{g+1} \to \Sigma_g$, then the map is instead given by
\begin{align}
  \rho_W: H^*(\Sym^{k+1}(\Sigma_{g+1})) & \to H^*(\Sym^k(\Sigma_g))\nonumber \\
  \omega & \mapsto \iota_c\omega,
\end{align}
where contraction is with respect to the intersection pairing on $H^1(\Sigma_g)$.

Having described Donaldson's ``TQFT", we now embark on proving formula (\ref{intersect}).

\section{Morse-Bott Gluing}

Here we state the appropriate Morse-Bott gluing results needed to obtain functoriality for our cobordisms.  Suppose we have (not necessarily elementary) cobordisms $W_0: \Sigma_0 \to \Sigma_1$ and $W_1: \Sigma_1 \to \Sigma_2$.  We have the asymptotic maps
\begin{align}
  \partial_\infty^0: M_{\eta_0}(W_0^*) & \to \V(\Sigma_0) \times \V(\Sigma_1)\\
  \partial_\infty^1: M_{\eta_1}(W_1^*) & \to \V(\Sigma_1) \times \V(\Sigma_2).
\end{align}
Define $\partial_\infty^{0,+}$ and $\partial_\infty^{0,-}$ to be $\partial_\infty^0$ and $\partial_\infty^1$  followed by projection onto the $\V(\Sigma_1)$ factor, respectively. The same analysis as in Lemma \ref{LemmaPert} shows that for generic compatible $\eta_0$ and $\eta_1$, i.e. those $\eta_0$ and $\eta_1$ satisfying our standard hypotheses from before, and which agree on the ends modeled on $\Sigma_1$, the map $\partial_\infty^0 \times \partial_\infty^1$ is transverse to the middle diagonal
$$\tilde \Delta := \V(\Sigma_0) \times \Delta_{\V(\Sigma_1) \times \V(\Sigma_1)} \times \V(\Sigma_2).$$
Thus, the preimage of the middle diagonal
\begin{equation}
  (\partial_\infty^0 \times \partial_\infty^1)^{-1}(\tilde \Delta) \label{eq:midD}
\end{equation}
is a smooth submanifold of $M_{\eta_0}(W_0^*) \times M_{\eta_1}(W_1^*)$.

\begin{Definition}\label{DefTrans}
  We say two cobordisms $W_0: \Sigma_0 \to \Sigma_1$ and $W_1: \Sigma_1 \to \Sigma_2$ are \textit{transverse} if the restriction maps $H^1(W_0) \to H^1(\Sigma_1)$ and $H^1(W_1) \to H^1(\Sigma_1)$ are transverse.
\end{Definition}

This definition is convenient because it implies that there are no gluing parameters when we glue the cobordisms $W_0$ and $W_1$, meaning that given $\spinc$ structures on $W_0$ and $W_1$ which are isomorphic when restricted to $\Sigma_1$, we get a unique glued up $\spinc$ structure on $W_1 \circ W_0$.  Our gluing results can be easily restated in the appropriate way for non-transverse cobordisms, by keeping track of the gluing parameter (much as in Theorem \ref{ThmTwo}), but for simplicity, we only state the results in the transverse case. 

\begin{Lemma}\label{Lemma-transglue}
Let $W_0: \Sigma_0 \to \Sigma_1$ and $W_1: \Sigma_1 \to \Sigma_2$ be two transverse cobordisms, where the boundary surfaces may be empty or disconnected. Let $W_T = W_0 \cup ([-T,T] \times \Sigma_1) \cup W_1$ be the composite cobordism with an added neck of length $2T$.  Then for $T$ sufficiently large, we have the following:
\begin{enumerate}
  \item The moduli space $M_{\eta_T}(W_T^*)$ is diffeomorphic to $(\partial_\infty^0 \times \partial_\infty^1)^{-1}(\tilde \Delta)$. Here, $\eta_T$ is the natural concatenation of the compatible perturbations $\eta_0$ and $\eta_1$ on $W_0$ and $W_1$.
  \item $[\partial_\infty M_{\eta_T}(W_T^*)]$ is homologous to $[\partial_\infty^1 M(W_1^*) \circ \partial_\infty^0 M(W_0^*)]$, where $\circ$ denotes geometric composition of the immersed Lagrangians.
\end{enumerate}
\end{Lemma}

\begin{proof}

(i) Our gluing construction in the Morse-Bott situation proceeds as follows.  Consider any cylindrical end cobordism $W^*$ and $u \in \fC^{s,\delta}(W^*)$. Smoothly identify  $T_{\partial_\infty(u)}\V(\Sigma)$, where $[0,\infty) \times \Sigma$ is the cylindrical end of $W^*$, with some space of configurations on $W^*$ that lie asymptotically in (\ref{eq2:kerB}) and which vanish outside of $[0,\infty) \times \Sigma$, i.e., extend elements of (\ref{eq2:kerB}) into $[0,\infty) \times \Sigma$ in some smooth way so that they extend into $W^*$ by being identically zero outside of the cylindrical end.  With this identification, define the space
\begin{equation}
  \tT^{s,\delta}_u = \tT^{s,\delta}(W^*) + T_{\partial_\infty(u)}\V(\Sigma). \label{eq:T+V}
\end{equation}
It is on these spaces where we can adapt the usual gluing methods in gauge theory to the appropriate operators.

Consider the map
\begin{align}
  \tT_{u_0}^{s,\delta}(W_0^*) \oplus \tT_{u_1}^{s,\delta}(W_1^*) & \to \tT^{s-1,\delta}(W_0^*) \oplus \tT^{s-1,\delta}(W_1^*) \oplus T_{a}\V(\Sigma_1) \nonumber \\
  (x_0,x_1) & \mapsto (\bH_{u_0}x_0, \bH_{u_1}x_1, (\partial_\infty^{0,+})_* x_0 - (\partial_\infty^{1,-})_* x_1) \label{op-MBglue}
\end{align}
where $u_i \in M_{\eta_i}(W_i^*)$, $i = 0,1$, are monopoles with matching limit $\a$.  The assumption that $M_{\eta_0}(W_0^*)$ and $M_{\eta_1}(W_1^*)$ are cut out transversally and that the maps $\partial_\infty^{0,+}: M(W_0^*) \to \V(\Sigma_1)$ and $\partial_\infty^{1,-}: M(W_1^*) \to \V(\Sigma_1)$ are transverse imply that (\ref{op-MBglue}) is surjective.  As we will show below, a right-inverse for this operator yields an approximate right-inverse for the operator
\begin{align}
    \bH_{u_0 \sharp_T u_1} : \tT^{s,\delta}_{u_0 \sharp_T u_1}(W_T^*) & \to \tT^{s;\delta}(W_T^*). \label{op-MBglue2}
\end{align}
Here the norm on the space $\tT^{s,\delta}_{u_0 \sharp_T u_1}(W_T^*)$ is defined in the usual way on the ends of $W_T^*$ but in addition has weight function $e^{\delta(T - |t|)}$ along the finite length neck $[-T,T]\times\Sigma_1$. This is the correct norm to use in our gluing construction since the function space on $W_T^*$ ought to approximate that on $W_0^* \coprod W_1^*$.  The configuration $u_0 \sharp_T u_1$ on $W_T^*$ is obtained by gluing $u_0$ and $u_1$ together as follows.  Let $\chi: [0,\infty) \to [0,1]$ be a smooth cutoff function which is identically one on $[0,1/4]$ and vanishes on $[1/2,\infty)$.  On the matching end $[0,\infty) \times \Sigma_1 \subset W_0^*$, define
$$u_0' = \chi_0(t/T)(u_0 - \g_\a) + \g_\a$$
where $\g_\a$ is the constant configuration equal to the limiting matching vortex $\a$ of $u_0$ and $u_1$, and $t \in [0,\infty)$.  Smoothly extend $u_0'$ to the rest of $W_i^*$ by setting it equal to $u_0$.  We define $\chi_1$ and $u_1'$ on $W_1^*$ similarly, which we can do since $W_1^*$ has the matching end $(-\infty,0] \times \Sigma_1$.  One can then concatenate the $u_i'$ in the obvious way by joining the $[0,T] \times \Sigma_1$ and $[-T,0] \times \Sigma_1$ portions of the matching ends, to obtain $u_0 \sharp_T u_1$ on $W_T$, which is identically $\a$ on the region $[-T/2, T/2] \times \Sigma_1$ inside the neck $[-T,T] \times \Sigma_1$.  To obtain an approximate right-inverse for (\ref{op-MBglue2}), which is the first step in the usual gluing story, one applies a similar splitting and gluing procedure to pass between (\ref{op-MBglue2}) and (\ref{op-MBglue}).  In detail:

Let $\varphi_0 + \varphi_1 = 1$ be a partition of unity on the interval $[-1,1]$, where $\varphi_0$ has support on $[-1,1/2]$ and $\varphi_1$ has support on $[-1/2,1]$. Then given $x \in \tT^{s;\delta}(W_T^*)$ we can split $x$ as $x = x_0 + x_1$, where $x_i = x$ on $W_i^*\setminus([0,\infty)\times\Sigma_1)$ and on the neck $[-T,T]\times\Sigma$, we have $x_i = \varphi_i(t/T)x$.
Including $(x_0,x_1)$ into $W_0^* \coprod W_1^*$ via extension by zero, one can then apply a right-inverse $\tilde R_T$ of (\ref{op-MBglue}) to this configuration to obtain $(\tilde x_0, \tilde x_1) := \tilde R_T(x_0,x_1,0)$.  Let $\tilde \varphi_0$ be a smooth function on $[0,2]$ with support on $[0,7/4]$ and which is identically one on $[0,3/2]$.  Define $\tilde \varphi_1$ on $[-2,0]$ by reflection.  Then truncating the $\tilde x_i$ along the necks $[0,\infty)\times \Sigma_1$ and $(-\infty,0]\times\Sigma$, respectively, via $(\tilde \varphi_i)(t/T)\tilde x_i$ and overlapping the two finite cylinders $[0,2T] \times \Sigma_1$ and $[-2T,0] \times \Sigma_1$ end-to-end so as to produce a cylinder of length $2T $ yields a glued together configuration $\tilde x_0 \tilde\sharp_T \tilde x_1$ on $W_T^*$. This defines for us an approximate right-inverse $Q_T$ for (\ref{op-MBglue2}):
$$Q_T(x) = \tilde x_0 \tilde\sharp_T \tilde x_1.$$
To check that we can perturb $Q_T$ to an honest right-inverse for the operator $\bH_{u_0\sharp_T u_1}$, we check that the operator norm of $\bH_{u_0\sharp_T u_1}Q_T - \mathrm{id}$ is small.  By construction, the $\tilde x_i$ decay exponentially along the corresponding necks of the $W_i^*$ to a common element in the tangent space to a vortex.  It follows that $\bH_{u_0\sharp_T u_1}Q_Tx$ is equal to $x = x_0 + x_1$ outside the neck $[-3T/4,3T/4] \times \Sigma_1$, and furthermore, on the neck, the difference between these two configurations has norm exponentially small in $T$ (if we set the norm of $x$ to be unity). Thus, the norm of $S:= \bH_{u_0\sharp_T u_1}Q_T - \mathrm{id}$ is exponentially small in $T$, and so we obtain a right-inverse $R_T := Q_T(1+S)^{-1}$ for $\bH_{u_0\sharp_T u_1}$.

Then the usual inverse function theorem methods allow us to use this right inverse to perturb the configuration $(u_0 \sharp_T u_1, 0)$ to a configuration $(u_0 \sharp_T u_1, 0) + (x,\xi)$, where $(x,\xi) \in \tT^{s,\delta}_{u_0\sharp_T u_1}$ solves
 \begin{equation}
   \begin{pmatrix}
     SW_3(u_0 \sharp_T u_1 + x) + \bd_{u_0 \sharp_T u_1 + x}\xi \\
      \bd_{u_0 \sharp_T u_1}^*x
   \end{pmatrix} = 0.
 \end{equation}
Indeed, this equation can be written as
$$\bH_{u_0\sharp_T u_1}(x,\xi) + q(x,\xi) = SW_3(u_0 \sharp_T u_1)$$
where $q$ is a quadratic multiplication map.  To solve this equation, it suffices to solve the equation
\begin{equation}
  y + q(R_T y) = SW_3(u_0 \sharp_T u_1). \label{eq:solve}
\end{equation}
If one traces through the construction, we have the upper bound $\|R_T\|_{\mathrm{Op}} \leq Ce^{\delta T}$. Here, we have an exponential growth estimate because an element of $T_{\partial_\infty(u)}\V(\Sigma)$ along the neck $[-T,T] \times \Sigma_1$ has norm $O(e^{\delta T})$ when regarded as an element of $\tT^{s,\delta}(W_T^*)$, whereas it has norm independent of $T$ regarded as an element belonging to the second factor of (\ref{eq:T+V}).  This difference in the way norms are defined on finite cylinders and on semi-infinite cylinders accounts for this exponential growth factor in tracing through the construction of $\tilde R_T$. Nevertheless, we can still apply the inverse function theorem because $SW_3(u_0 \sharp_T u_1) = O(e^{-\mu T})$ decays rapidly, where we can choose any $\mu$ such that $\delta < \mu < \delta_0$, where $\delta_0$ is the constant appearing in Lemma \ref{LemmaExpDecay}.  Indeed, recall that $\delta > 0$ is a sufficiently small constant fixed at the outset. Then applying the inverse function theorem with estimate tells us that for sufficiently large $T$, there is a unique solution $y$ to (\ref{eq:solve}) with $y = O(e^{-\mu T})$. Such a solution automatically satisfies $\xi = 0$, since the terms $SW(u_0 \sharp_T u_1+x)$ and $\bd_{u_0 \sharp_T u_1 + x}\xi$ live in complementary subbundles of $\tT^{s-1,\delta}$ (cf. Lemma \ref{LemmaTw}).  In this way, we see that $u_0 \sharp_T u_1 + x$ is a monopole on $W_T$.  The same arguments as in \cite[Chapter 4.4]{Don} show that this construction works in families, so that a family of pairs of monopoles $(u_0(s), u_1(s))$ with matching limits can be perturbed to yield a family of monopoles on $W_T$, this correspondence being a smooth bijection. Moreover, the surjectivity of this gluing construction also follows from the same arguments as in \cite{Don}.  This establishes the claimed diffeomorphism.  

(ii) Each large $T$ gives rise to a diffeomorphism from the fiber product (\ref{eq:midD}) to $M(W_T^*)$.  Composing this diffeomorphism with $\partial_\infty: M(W_T^*) \to \V(\Sigma_0 \times \Sigma_2)$, the composite map varies smoothly with $T$ and converges to the map $\partial_\infty^0 \times \partial^1_\infty$ defined on the fiber product.  This provides the required homological equivalence.
\end{proof}

The above lemma immediately implies the composition rule:

\begin{Corollary}
  For $W_0$ and $W_1$ as above, we have
  $$\varrho_{W_0 \cup W_1} = \varrho_{W_1} \circ \varrho_{W_0}.$$
\end{Corollary}

\begin{proof}
  We have
\begin{align*}
  \varrho_{W_1} \circ \varrho_{W_0} &= (\partial_\infty^1)_*[M(W_1^*)] \circ (\partial_\infty^0)_*[M(W_0^*)]\\
  &= [\partial_\infty^1(M(W_1^*)) \circ \partial_\infty^0(M(W_0^*))]
\end{align*}
where the last line denotes the map in $\Hom(H_*(\V(\Sigma_0)), H_*(\V(\Sigma_1))$ induced from the homology class of the immersed Lagrangian submanifold $\partial_\infty^1(M(W_1^*)) \circ \partial_\infty^0(M(W_0^*)$ of $-\V(\Sigma_0) \times \V(\Sigma_2)$.  Indeed, the last equality is an exercise in the intersection pairing on homology and follows straight from the definitions.  On the other hand, $\partial_\infty^1(M(W_1^*)) \circ \partial_\infty^0(M(W_0^*)$ is homologous to $\partial_\infty(W_T)$ for large $T$, and the latter induces the map $\rho_{W_0 \cup W_1}$.
\end{proof}

\section{Proof of Theorem \ref{ThmTwo}}

Stretch the metric along $\Sigma$ so that we may consider the closed manifold $$Y_T = (Y\setminus \Sigma) \cup ([-T,T] \times \Sigma) =: W_0 \cup W_{1,T}$$ obtained by replacing $\Sigma$ with the tube $[-T,T]\times\Sigma$.  From the Morse-Bott gluing lemma, for large $T$, the space of monopoles on $Y_T$ is in bijective correspondence with the fiber product of
\begin{align*}
  \partial_\infty^0: M(W_0^*) & \to \V(-\Sigma \coprod \Sigma)\\
  \partial_\infty^1: M(W_1^*) & \to \V(-\Sigma \coprod \Sigma), \qquad W_1 := W_{1,1} = [-1,1]\times\Sigma
\end{align*}
and hence to a discrete set of points after arranging for a transverse intersection of the resulting immersed Lagrangians (where as usual, suitable perturbations on $Y_T$ have been chosen, in particular, those which are product when restricted to the neck $[-T,T] \times \Sigma$ so that they extend in the obvious way when cylindrical ends are attached to $W_0$ and $W_{1,T}$).  As explained previously, because we are considering the intersection in $\V(-\Sigma \coprod \Sigma)$, the monopoles we obtain on $Y_T$ are with respect to the set of $\spinc$ structures obtained from all possible gluing parameters.
At a monopole $u$ on $Y_T$, the extended Hessian $\bH_u$ is invertible, and the contribution of $u$ to the Seiberg-Witten invariant on $Y_T$ is obtained by comparing the canonical orientation of $\det \bH_u$ with that given by the homology orientation on $Y_T$.  Up to an overall sign in the Seiberg-Witten invariant (as determined by the homology orientation), we will show that this comparison, i.e. the relative sign difference, is precisely given by the signed intersection of the oriented Lagrangian subspaces $\partial_\infty^0 T_{u_0}M(W_0^*)$ and $\partial_\infty^1 T_{u_1}M(\R \times \Sigma)$ inside $T_{\partial_\infty^0 u_0}\V(-\Sigma  \coprod \Sigma)$. Call this latter sign $\eps(u_0,u_1)$.  We prove the following key lemma:

\begin{Lemma}\label{LemmaLineMult}
  Let $U_i$ be precompact open subsets  of $\fC^{s;\delta}(W_i^*)$.
\begin{enumerate}
  \item For sufficiently large $T$, there exists a continuous isomorphism of determinant lines
\begin{equation}
   \det (\bH_{u_0}) \otimes \det (\bH_{u_1}) \to \det (\bH_{u_0 \sharp_T u_1}), \label{linemult}
\end{equation}
where $u_i \in U_i$ are such that $\partial_\infty^0u_0 = \partial_\infty^1u_1$.
  \item Choose $T$ sufficiently large so that for the $u_i$ being any pair of monopoles with matching limits, the operator $\bH_{u_0 \sharp_T u_1}$ is invertible. Then by a suitable choice of orientations for the determinant lines for the $\fC^{s;\delta}(W_i^*)$, the induced orientation on $\det (\bH_{u_0 \sharp_T u_1})$ coming from (\ref{linemult}) differs from the canonical orientation of $\det (\bH_{u_0 \sharp_T u_1})$ by $\eps(u_0,u_1)$.
\end{enumerate}
\end{Lemma}
Here, we suppress from our notation the glued $\spinc$ structure on $Y_T$ we obtain from the matching pair $(u_0,u_1)$.

Before proving the lemma, let's see how the lemma proves the theorem.  The continuity of (\ref{linemult}) gives us a trivialization of the determinant line of $\fC(Y_T)$ over the set of all points of the form $u_0 \sharp_T u_1$.  One can arrange the $U_i$ so that they contain every monopole on $W_i^*$ up to gauge, and furthermore, that the associated set of elements $U_0 \sharp_T U_1 := \{u_0 \sharp_T u_1\}$ is a connected subset of $\fC(Y_T)$. Thus, a trivialization of the determinant line over $U_0 \sharp_T U_1$ induces an orientation of the determinant line on all of $\fC(Y_T)$, from which passing to the quotient, we get an induced trivialization of the determinant line over all of $\fB(Y_T)$ (since the determinant line is globally a trivial line bundle).  For large $T$, given two monopoles $u_i$, one can join $u_0 \sharp_T u_1$ to its associated monopole $u$ under gluing by a short path of configurations on $Y_T$, all of whose extended Hessian operators are invertible. (That this is possible follows from the proof of Lemma \ref{Lemma-transglue}, which shows that the distance between the monopole $u$ and the approximate monopole $u_0 \sharp_T u_1$ is much smaller than the operator norm of the inverse of $\bH_{u_0\sharp_T u_1}$).  Thus, our trivialization of $\fB(Y_T)$, induced by (\ref{linemult}), differs from the canonical orientation at each monopole associated to $u_0\sharp_T u_1$ by the sign $(-1)^{\eps(u_0,u_1)}$, which therefore gives us the correct signed count in the Seiberg-Witten invariant if we homology orient the determinant line on $\fB(Y_T)$ in the corresponding way.  Note that in fixing a homology orientation for $Y$, we obtain an orientation for $\fB(Y_T)$ with respect to any $\spinc$ structure in the usual way, these orientations all being compatible with respect to different gluing parameters since the gluing parameters act as complex linear maps on spinors.  Thus, the preceding analysis shows us that for the appropriate choice of homology orientation on $Y$, the signed intersection of the Lagrangians $\partial_\infty^0M(W_0^*)$ and $\partial_\infty^1M(W_1^*)$ agrees exactly with the signed count of monopoles on $Y_T$ arising from Seiberg-Witten invariant summed over all possible glued $\spinc$ structures.

Thus, to establish the main theorem, it remains to establish Lemma \ref{LemmaLineMult}.\\

\noindent \textsc{Proof of Lemma \ref{LemmaLineMult}}: (i) Consider the map
\begin{align}
  \tT_{u_0}^{s,\delta}(W_0^*) \oplus \tT_{v_1}^{s,\delta}(W_1^*) & \to \tT^{s-1,\delta}(W_0^*) \oplus \tT^{s-1,\delta}(W_1^*) \oplus T_{\a}\V(-\Sigma \coprod \Sigma) \nonumber \\
  (x_0,x_1) & \mapsto (\bH_{u_0}x_0, \bH_{u_1}x_1, (\partial_\infty^0)_* x_0 - (\partial_\infty^1)_* x_1) \label{op-partial}
\end{align}
as in (\ref{op-MBglue}), where $\a = \partial_\infty^0 u_0 = \partial_\infty^1 u_1$. Denote the operator (\ref{op-partial}) by $L_{u_0,u_1}$. We will show that there is an isomorphism from the determinant line of $L_{u_0,u_1}$ to $\det \bH_{u_0 \sharp_T u_1}$ on $Y_T$.  The multiplication (\ref{linemult}) will then follow since the determinant line of (\ref{op-partial}) is isomorphic to $\det (\bH_{u_0}) \otimes \det (\bH_{u_1})$, by considering the linear homotopy $t(\partial_\infty^0 - \partial_\infty^1)$, $t \in [0,1]$ on the third factor. Here, we use that the top exterior power of the vortex tangent space is canonically oriented since it is a complex vector space.

So suppose we are given any $u_i \in U_i$ with matching limits. Then we can find a finite dimensional space of smooth configurations $Z$ compactly supported in the interior of $W_0 \coprod ([-1,1] \times \Sigma)$ such that if
$$S: Z \to \tT^{s-1,\delta}(W_0^*) \oplus \tT^{s-1,\delta}(W_1^*),$$
denotes the inclusion map, then both $L_{u_0,u_1} + S$ and $\bH_{u_0\sharp_T u_1} + S$ are surjective. (The domains of both these operators have been enlarged to contain $Z$, and in the latter case, $Z$ also sits naturally inside the range of $\bH_{u_0 \sharp_T u_1}$ since the $Z$ are supported on $W_0 \coprod W_1 \subset W_T = Y_T$. Hence, by adding $S$, we have increased the index by $\dim Z$).  Indeed, observe from unique continuation and the Fredholm property of $L_{u_0,u_1}$ that the $L^2$-orthogonal complement of the range of $L_{u_0,u_1}$ is a finite dimensional subspace of $\tT^{s-1,\delta}(W_0^*) \oplus \tT^{s-1,\delta}(W_1^*)$ consisting of configurations that do not vanish identically on any open set.  It follows that the choice of $S$ as above making $L_{u_0,u_1} + S$ surjective is possible, and furthermore, this $S$ works for all $u_i$ nearby.  To see that this same $S$ also makes $\bH_{u_0\sharp_T u_1} + S$ surjective, for sufficiently large $T$, one simply repeats the same arguments as in the Morse-Bott gluing construction of the previous section to see that a right inverse for $L_{u_0,u_1} + S$ gives a right-inverse for $\bH_{u_0\sharp_T u_1} + S$.

It suffices now to construct an isomorphism between the determinant lines for the stabilized maps $\bH_{u_0,u_1} + S$ and $\L_{u_0,u_1} + S$ for all $u_i$ on a sufficiently small open set for which the stabilization map $S$ provides us with surjective operators.  In fact, under these conditions, we construct an isomorphism between the kernels of the stabilized operators. Indeed, any element of the kernel of $L_{u_0,u_1} + S$ can be glued together (since their asymptotic limits agree) to give an element approximately in the kernel of $\bH_{u_0 \sharp_T u_1}$, and then one uses the surjectivity of $\bH_{u_0 \sharp_T u_1}$ to perturb to an exact solution.  In detail, using the notation used within the proof of Lemma \ref{Lemma-transglue}, given $z \in Z$ and $\tilde x_i \in \tT^{s,\delta}(W_i^*) + T_{\partial_\infty(u_i)}\V$ such that $(\tilde x_0, \tilde x_1) + z \in \ker (\L_{(u_0,u_1)} + S)$, we can define an element $x$ belonging to $\tT^{s,\delta}(W_T)$ such that
\begin{align*}
  x_0 := x|_{W_0 \cup \left([-T,0] \times (\Sigma \coprod \Sigma)\right)} = \chi_0(t/T)(\tilde x_0 - v_\infty) + v_\infty \\
  x_1 := x|_{[0,T] \times (\Sigma \coprod \Sigma)} = \chi_1(t/T)(\tilde x_1 - v_\infty) + v_\infty
\end{align*}
where $v_\infty$ is the $T_{\partial_\infty(u_i)}\V$ component, i.e., the linearized vortex component of the $\tilde x_i$. From the way our norms are defined, we have the bounds
\begin{equation}
  C\|x\|_{W_0^* \coprod W_1^*} \leq \|(\tilde x_0, \tilde x_1)\|_{\tT^{s,\delta}(W_T)} \leq Ce^{\delta T}\|x\|_{W_0^* \coprod W_1^*} \label{eq:ker1}
\end{equation}
where $C$ is a constant independent of $T$ for $T$ large.  Here, the middle norm in the above denotes the norm on $\oplus_i (\tT^{s,\delta}(W_i^*) \oplus T_{\partial_\infty(u_i)}\V$. (We have chosen some fixed but otherwise arbitrary norm on the finite dimensional space $T_{\partial_\infty(u_i)}\V$.) The first inequality comes from the fact that elements of the kernel of $(\L_{(u_0,u_1)} + S)$ decay exponentially along the end.

On the other hand, proceeding as in the proof of Lemma \ref{Lemma-transglue}, we have the error estimate
\begin{equation}
  \|\bH_{u_0\sharp_Tu_1}x + z\|_{\tT^{s,\delta}(W_T)} \leq Ce^{(-\mu + \delta)T}\|(\tilde x_0, \tilde x_1)\|_{W_0^* \coprod W_1^*} \label{eq:ker2}
\end{equation}
for some $\mu > \delta$ to be determined. Indeed, it suffices to examine $\bH_{u_0\sharp_Tu_1}x + z$ on $W_0 \cup \left([-T,0] \times (\Sigma \coprod \Sigma)\right)$ (the estimate on the other half of $Y_T$ is similar). We have $\bH_{u_0\sharp_Tu_1} = \bH_{u_0} + O(e^{-\mu T})$ since $u_0$ converges exponentially fast to a vortex, where $\mu < \delta_0$ as in Lemma \ref{LemmaExpDecay}.  Thus it suffices to estimate
\begin{align}
  \bH_{u_0}x + z &= \bH_{u_0}x_0 + z \nonumber \\
  &= \left(\bH_{u_0}\chi_0(t/T)\right)(\tilde x_0 - v_\infty) + \chi_0(t/T)\left(\bH_{u_0}(\tilde x_0 - v_\infty) + \bH_{u_0}v_\infty + z\right) \nonumber \\
  & \qquad + (1 - \chi_0(t/T))\bH_{u_0}v_\infty. \label{eq:expsmall}
\end{align}
The first term of (\ref{eq:expsmall}) is exponentially small since $\chi_0(t/T)$ has support on $[-T/4,T]$ and $(\tilde x_0 - v_\infty)$ decays exponentially fast.  In fact, since $L_{(u_0,u_1)}(\tilde x_0, \tilde x_1) = 0$ on the cylindrical ends of $W_0^* \coprod W_1^*$, $\tilde x_0 - v_\infty$ decays at any exponential rate less than the spectral gap of $\mathsf{B}_\a$ (i.e. the distance from $0$ to the first non-zero eigenvalue), where $\a$ is the limiting vortex of the $u_i$. Since the vortex moduli space is compact, this spectral gap is uniformly bounded away from zero in $\a$, and so when $\delta > 0$ is sufficiently small, the first term of (\ref{eq:expsmall}) is compatible with the estimate (\ref{eq:ker2}) for some $\mu > \delta$.  Next, the second term of (\ref{eq:expsmall}) is identically zero since $(\tilde x_0,\tilde x_1) + z \in \ker (L_{u_0,u_1}+S)$.  Finally the third-term is exponentially small since
$$(1 - \chi_0(t/T))\bH_{u_0}v_\infty = (1 - \chi_0(t/T))(u_0 - \g_\a)v_\infty$$
and $u_0 - \g_\a$ is $O(e^{-\mu t})$, where $\mu < \delta_0$. (Note that $\delta_0$ is precisely the infimum of the spectral gap discussed above.) Altogether, having estimated the three terms of (\ref{eq:expsmall}), this establishes (\ref{eq:ker2}).

It follows from (\ref{eq:ker1}) and (\ref{eq:ker2}), and from the operator norm of $(\bH_{u_0\sharp_Tu_1} + S)^{-1}$ being bounded by $O(e^{\delta T})$, that the map
\begin{align*}
  \ker (\L_{(u_0,u_1)} + S) & \to \ker (\bH_{u_0\sharp_Tu_1}+S) \\
  (\tilde x_0,\tilde x_1) + z & \mapsto x + z - (\bH_{u_0\sharp_Tu_1}+S)^{-1}(\bH_{u_0\sharp_Tu_1}+S)(x + z)
\end{align*}
is injective for sufficiently large $T$ and $\mu > 2\delta$. Since $\L_{(u_0,u_1)} + S$ and $\bH_{u_0\sharp_Tu_1}+S$ both have the same index and are surjective, by construction, it follows that the above map of kernels is an isomorphism. Passing to determinants, this establishes the requisite multiplication map of determinant lines for $(u_0,u_1)$ belonging to small open subsets of the $U_i$.  Since the $U_i$ are precompact, we can choose a single fixed large enough $T$ so that we have the requisite multiplication on open subsets of $U_i$ belonging to some finite cover of the $U_i$.  This gives us the multiplication map on all of the $U_i$, well defined at the level of orientations (i.e. well-defined modulo an overall positive scaling of the determinant lines) since all the different choices involved in the above multiplication (choice of cut-off functions, stabilization map, etc) are all homotopic. (That the multiplication is well-defined up to orientation is enough for our purposes. However, one can then patch together these orientation-compatible multiplications on the open cover of the $U_i$ to yield a well-defined multiplication map which is fiberwise an isomorphism over the $U_i$ although this final step is not necessary).

(ii) It is now a matter to see how this gluing construction behaves when we glue the determinant line of two monopoles with matching limits.  Observe that in this case, the operator (\ref{op-partial}) is an isomorphism by hypothesis, since the monopoles are cut out transversally and the tangent spaces to the monopole moduli spaces are transverse Lagrangians at infinity.  Thus, the determinant line of (\ref{op-partial}) carries a canonical orientation, which induces the canonical orientation of $\det (\bH_{u_0\sharp_Tu_1})$ when we glue.  On the other hand, consider the determinant line along the path of operators $P_t := (\bH_{u_0}x_0, \bH_{u_1}x_1, t(\partial_\infty^0 x_0 - \partial_\infty^1 x_1))$, $t \in [0,1]$.  At $t=0$, an orientation for $P_0$ is determined by an orientation for $\ker \bH_{u_0} \oplus \ker \bH_{u_1}$, since the tangent space to the space of vortices is a complex vector space and hence canonically oriented.  Orienting each $\ker \bH_{u_i}$ (which is the same as orienting the moduli spaces $\M(W_i^*)$) to obtain an orientation of $\det P_0$, orientation transport along $[0,1]$ shows that the induced orientation on $\det P_1$ differs from the canonical orientation by $\partial_\infty^0(\ker \bH_{u_0}) \cap \partial_\infty^1(\ker \bH_{u_1})$, that is, the signed intersection of the corresponding oriented Lagrangian subspaces of $T_{[\a]}\V$.  This establishes the claim about signs when we glue.\End





\begin{thebibliography}{BBLZ}

\bibitem{APS} M. F. Atiyah, V. K. Patodi, and I. M. Singer, I. M. \textit{Spectral asymmetry and Riemannian geometry. I.} Math. Proc. Cambridge Philos. Soc. 77 (1975), 43–-69.
\bibitem{BBLZ} B. Booss-Bavnbek, M. Lesch, and C. Zhu. \textit{The Calder\'on Projection: New Definition and Applications.}  J. Geom. Phys.  59  (2009),  no. 7, 784--826.
\bibitem{Don} S. Donaldson. \textit{Floer homology groups in Yang-Mills theory.} With the assistance of M. Furuta and D. Kotschick. Cambridge Tracts in Mathematics, 147. Cambridge University Press, Cambridge, 2002.
\bibitem{Don99} S. Donaldson. \textit{Topological field theories and formulae of Casson and Meng-Taubes.} Proceedings of the Kirbyfest (Berkeley, CA, 1998), 87–-102, Geom. Topol. Monogr., 2, Geom. Topol. Publ., Coventry, 1999.
\bibitem{DK} S. Donaldson and P. Kronheimer. \textit{The geometry of four-manifolds.}
    Oxford Mathematical Monographs. The Clarendon Press, Oxford University Press, New York, 1990.
\bibitem{HL} M. Hutchings and Y.-J. Lee. \textit{Circle-valued Morse theory and Reidemeister torsion.}
    Geom. Topol. 3 (1999), 369–-396.
\bibitem{GP} O. Garcia-Prada. \textit{A direct existence proof for the vortex equations over a compact Riemann surface.} Bull. London Math. Soc. 26 (1994), no. 1, 88–-96.
\bibitem{KM} P. Kronheimer and T. Mrowka. \textit{Monopoles and 3-manifolds.} Cambridge University Press, Cambridge, 2007.
\bibitem{KLT} C. Kutluhan, Y.-J. Lee, and C. Taubes. \textit{HF=HM I : Heegaard Floer homology and Seiberg--Witten Floer homology.} arxiv:1007.1979.
\bibitem{Lip} M. Lipyanskiy. \textit{A Semi-Infinite Cycle Construction of Floer Homology.} PhD Thesis, 2008.
\bibitem{LM} R. Lockhart and R. McOwen. \textit{Elliptic differential operators on noncompact manifolds.}
    Ann. Scuola Norm. Sup. Pisa Cl. Sci. (4) 12 (1985), no. 3, 409–-447.
\bibitem{Mark} T. Mark. \textit{Torsion, TQFT, and Seiberg-Witten invariants of 3-manifolds.} Geom. Topol. 6 (2002), 27–-58.
\bibitem{MT} G. Meng and C. Taubes. \textit{SW = Milnor Torsion.} Math. Res. Lett. 3 (1996), 661–-674.
\bibitem{MMR} J. Morgan, T. Mrowka, D. Ruberman. \textit{The $L^2$-moduli space and a vanishing theorem for Donaldson polynomial invariants.} Monographs in Geometry and Topology, II. International Press, Cambridge, MA, 1994.
\bibitem{MMS} J. Morgan, T. Mrowka, Z. Szabo. \textit{Product formulas along $T^3$ for Seiberg-Witten invariants.} Math. Res. Lett. 4 (1997), no. 6, 915-–929.
\bibitem{MST} J. Morgan, Z. Sz\'{a}bo, and C. Taubes. \textit{A product formula for the Seiberg-Witten invariants and the generalized Thom conjecture.} J. Differential Geom. 44 (1996), no. 4, 706–-788.
\bibitem{MOY} T. Mrowka, P. Ozsvath, B. Yu. \textit{Seiberg-Witten monopoles on Seifert fibered spaces.}
    Comm. Anal. Geom. 5 (1997), no. 4, 685–-791.
\bibitem{N} T. Nguyen. \textit{The Seiberg-Witten Equations on Manifolds with Boundary.} PhD Thesis, 2011.
\bibitem{N1} T. Nguyen. \textit{The Seiberg-Witten Equations on Manifolds with Boundary I.} Comm. Anal. Geom. 20 (2012), no. 3, 565--676.
\bibitem{N2} T. Nguyen. \textit{The Seiberg-Witten Equations on Manifolds with Boundary II.} arxiv: 1008.2017
\bibitem{Nic} L. Nicolaescu. \textit{Notes on Seiberg-Witten theory}. Graduate Studies in Mathematics, 28. American Mathematical Society, Providence, Rhode Island, 2000.
\bibitem{Nic95} L. Nicolaescu. \textit{The Maslov index, the spectral flow, and decompositions of manifolds}. Duke Math. J. 80 (1995), no. 2, 485-–533.
\bibitem{Sal} D. Salamon. \textit{Seiberg-Witten invariants of mapping tori, symplectic fixed points, and Lefschetz numbers.} Proceedings of 6th G\"okova Geometry-Topology Conference. Turkish J. Math. 23 (1999), no. 1, 117-–143.
\bibitem{T} C. Taubes. \textit{Casson's invariant and gauge theory.} J. Differential Geom. 31 (1990), no. 2, 547-–599.
\bibitem{Tur} V. Turaev. \textit{A combinatorial formulation for the Seiberg-Witten invariants of 3-manifolds.}
Math. Res. Lett. 5 (1998), no. 5, 583–-598.

\end{thebibliography}
\end{document}